\documentclass[a4paper,12pt,reqno]{amsart}
\usepackage[T1]{fontenc}
\usepackage[utf8]{inputenc}
\usepackage[english]{babel}
\usepackage[left=25mm,right=25mm,top=35mm,bottom=35mm]{geometry}
\usepackage{amsmath,amssymb,amsthm}
\usepackage{pgfplots}
\usepackage{pgfplotstable}
\usepackage[margin=0pt,font={small}]{caption}
\usepackage{booktabs}
\usepackage{color,colortbl}
\usepackage[colorlinks=true,allcolors=blue,breaklinks]{hyperref}

\newcommand{\y}{\mathbf{y}} 
\newcommand{\z}{\mathbf{z}} 

\newcommand{\G}{\Gamma}
\newcommand{\N}{\mathbb{N}}

\newcommand{\R}{\mathbb{R}}
\newcommand{\V}{\mathbb{V}}
\newcommand{\X}{\mathbb{X}}
\newcommand{\Y}{\mathbb{Y}}

\newcommand{\MM}{\mathcal{M}}
\newcommand{\NN}{\mathcal{N}}

\renewcommand{\SS}{\mathcal{S}}
\newcommand{\TT}{\mathcal{T}}

\DeclareMathOperator*{\refine}{refine}

\newcommand{\norm}[3][]{#1\|#2#1\|_{#3}}

\newcommand{\dx}{\mathrm{d}x}

\newcommand{\Colpts}{\mathcal{Y}} 
\newcommand{\mf}{\kappa} 
\newcommand{\indset}{\Lambda} 
\newcommand{\markindset}{\Upsilon} 
\newcommand{\rmarg}{{\rm R}} 
\newcommand{\nnu}{\boldsymbol{\nu}} 
\newcommand\1{\boldsymbol{1}}
\newcommand{\LagrBasis}[2]{L_{#1}^{#2}} 
\newcommand{\LagrBasisHat}[2]{\widehat L_{#1}^{#2}} 
\newcommand{\scsol}{u_{\bullet}^{\rm SC}} 

\pgfplotsset{
every axis/.append style={
font={\fontsize{8pt}{12pt}\selectfont},  
},
tick label style={font=\tiny},
title style={font=\tiny,yshift=-1.5ex},
xlabel style={font=\tiny,yshift=+1.0ex},
ylabel style={font=\tiny,yshift=-1.2ex},
}

\definecolor{myBrown}{rgb}{0.6 0.4 0.2}
\definecolor{myOrange}{rgb}{1.0 0.6 0.2}
\definecolor{myLightGray}{RGB}{235,235,235}
\definecolor{myViolet}{RGB}{153,50,204}
\newtheorem{theorem}{Theorem}

\newtheorem{algorithm}[theorem]{Algorithm}

\makeatletter
\def\@seccntformat#1{%
  \protect\textup{\protect\@secnumfont
    \ifnum\pdfstrcmp{subsection}{#1}=0 \bfseries\fi
    \csname the#1\endcsname
    \protect\@secnumpunct
  }%
}
\makeatother
\newcommand*\patchAmsMathEnvironmentForLineno[1]{%
  \expandafter\let\csname old#1\expandafter\endcsname\csname #1\endcsname
  \expandafter\let\csname oldend#1\expandafter\endcsname\csname end#1\endcsname
  \renewenvironment{#1}%
     {\linenomath\csname old#1\endcsname}%
     {\csname oldend#1\endcsname\endlinenomath}}%
\newcommand*\patchBothAmsMathEnvironmentsForLineno[1]{%
  \patchAmsMathEnvironmentForLineno{#1}%
  \patchAmsMathEnvironmentForLineno{#1*}}%
\AtBeginDocument{%
\patchBothAmsMathEnvironmentsForLineno{equation}%
\patchBothAmsMathEnvironmentsForLineno{align}%
\patchBothAmsMathEnvironmentsForLineno{flalign}%
\patchBothAmsMathEnvironmentsForLineno{alignat}%
\patchBothAmsMathEnvironmentsForLineno{gather}%
\patchBothAmsMathEnvironmentsForLineno{multline}%
}
\usepackage[mathlines]{lineno}

\usepackage{fancyhdr}
\advance\footskip0.5cm
\definecolor{otherblue}{rgb}{0,0.3,0.6}
\def\rbl#1{\textcolor{black}{#1}}
\newcommand\abrev[1]{{\color{black}#1}}

\newcommand\partI{{\textcolor{blue}{\cite{bsx21}}}}
\newcommand\abrevx[1]{{\color{black}#1}}
\newcommand\rblx[1]{{\color{black}#1}}

\title{Error estimation and adaptivity for stochastic collocation finite elements\\ Part II: 
\rblx{multilevel} approximation}
\author{Alex Bespalov}
\address{School of Mathematics, University of Birmingham, Edgbaston, Birmingham B15 2TT, UK}
\email{a.bespalov@bham.ac.uk}
\author{David J. Silvester}
\address{Department of Mathematics, University of Manchester, Oxford Road, Manchester M13 9PL, UK}
\email{d.silvester@manchester.ac.uk}
\date{\today}

\begin{document}

\begin{abstract}
A  multilevel adaptive refinement strategy for solving linear elliptic partial differential \rblx{equations} with random data is  
recalled in this work. The strategy extends the a posteriori error estimation framework introduced by  Guignard \& Nobile in 2018 ({\sl SIAM J. Numer. Anal.}, {\bf 56}, 3121--3143)  
to cover  problems with  a {\it nonaffine} parametric coefficient dependence. A suboptimal, but nonetheless reliable and convenient implementation of the strategy involves approximation of the decoupled  PDE problems with a common finite element approximation space.  Computational results obtained using  such a  {\it single-level} strategy are presented in 
part I of this work (Bespalov, Silvester \& Xu, \href{http://arXiv.org/abs/2109.07320}{\hbox{arXiv:2109.07320}}). Results obtained using a potentially more efficient {\it multilevel} approximation strategy, where meshes are individually tailored,  are discussed herein.
 The codes used to generate the numerical results are available online. 
\end{abstract}

\maketitle
\thispagestyle{fancy}

\section{Introduction} \label{sec:intro}

Partial differential equations (PDEs) with uncertain inputs have provided engineers and scientists with enhanced 
fidelity  in the modelling of real-life phenomena, especially within the last decade.
Sparse grid stochastic collocation representations of parametric uncertainty, in combination with finite element 
discretization of physical space, have emerged as an efficient alternative to Monte-Carlo strategies over this period, 
especially in the context of nonlinear PDE models or linear PDE problems that are nonlinear in the parameterization
of the uncertainty.

The combination of adaptive sparse grid methods with a  hierarchy of spatial approximations is a relatively new 
development, see for example,~\cite{LangSS20,TeckentrupJantschWebsterGunzburger2012}.   In our precursor 
paper~\cite{bsx21} (part I), we extended  the  adaptive framework developed by  Guignard \& Nobile~\cite{GuignardN18} 
and  presented a critical comparison of alternative strategies in the context of solving a model 
problem that combines strong anisotropy in the parametric dependence with singular behavior in the physical space.
The numerical results presented in \partI\ demonstrate the effectivity and  robustness
of our error estimation strategy as well as the utility of the error indicators guiding the adaptive refinement process.
The results  in \partI\ also showed that optimality of convergence is difficult to achieve using a simple single-level approach 
where a single finite element space is associated with all active collocation points.
The main aim of this contribution is to see if optimal convergence rates can be recovered by computing results using
a multilevel implementation of the  algorithm  outlined in \partI.

The convergence of a modified version of the adaptive algorithm in~\cite{GuignardN18} has been established by 
Eigel et al.~\cite{eest20} and independently by Feischl \& Scaglioni~\cite{FeischlS21}. The authors of~\cite{FeischlS21} 
note that the main difficulty in establishing convergence  is ``the interplay of parametric refinement and finite element 
refinement''. This interplay is the focus of this contribution.

The model problems that are of interest  are stated in section~\ref{sec:problem}. The only difference
from  the problem statement in \partI\  is that we also  cover the case where the right-hand side function
\abrevx{has a parametric dependence}.
The adaptive solution  algorithm from \partI\ is extended to cover the case of a 
non-deterministic right-hand side function  in section~\ref{sec:scfem}. The novel contribution of this work 
 primarily  lies in section~\ref{sec:results}, where we compare  numerical results  obtained with our multilevel 
 algorithm with those generated using  a single-level strategy
 \abrevx{and with those computed using the multilevel stochastic Galerkin finite element method (SGFEM)}.

\section{\rbl{Parametric model problems}} \label{sec:problem}

Let $D \in \R^2$ be a bounded Lipschitz domain with polygonal boundary $\partial D$.
Let $\Gamma := \G_1 \times \G_2 \times \cdots \times \G_M$ denote the parameter domain in $\R^M$,
where $M \in \N$ and each $\G_m$ ($m = 1,\ldots,M$) is a bounded interval in~$\R$. 
We introduce a probability measure $\pi(\y) := \prod_{m=1}^M \pi_m(y_m)$ on $(\G,\mathcal{B}(\G))$;
here, $\pi_m$ denotes a 
Borel probability measure on $\G_m$ ($m = 1,\ldots,M$) and
$\mathcal{B}(\G)$ is the Borel $\sigma$-algebra on $\G$.


The first model  problem is the parametric elliptic problem analyzed in \partI:
we seek $u \colon \overline D \times \G \to \R$ satisfying
\begin{subequations}
\begin{align} 
\label{eq:pde:strong}
\begin{aligned}
-\nabla \cdot (a(\cdot, \y)\nabla u(\cdot, \y))
&= f 
&& \text{in $	D$},\\ 
u(\cdot, \y) &= 0  && \text{on $\partial D$,} 
\end{aligned}
\\
\intertext{$\pi$-almost everywhere on $\G$. The second  model problem is
to find $u \colon \overline D \times \G \to \R$ satisfying}
\label{eq:pde2:strong}
\begin{aligned}
-\nabla^2 u(\cdot, \y))
&= f({\cdot,\y})
&& \text{in $	D$},\\ 
u(\cdot, \y) &= 0  && \text{on $\partial D$,} 
\end{aligned}
\end{align}
\end{subequations}
$\pi$-almost everywhere on $\G$.

\rblx{In} the first model problem, the deterministic right-hand side function $f \in L^2(D)$ and
the coefficient $a$ is a random field on $(\G,\mathcal{B}(\G),\pi)$ over $L^\infty(D)$.
In this case we will assume that there exist constants $a_{\min},\, a_{\max}$ such that
\begin{equation} \label{eq:amin:amax}
   0 < a_{\min} \le \operatorname*{ess\;inf}_{x \in D} a(x,\y) \le \operatorname*{ess\;sup}_{x \in D} a(x,\y) \le a_{\max} < \infty \quad
   \text{$\pi$-a.e. on $\G$}.
\end{equation}
This assumption  implies the following norm equivalence:
for any $v \in \X := H^1_0(D)$ there holds
\begin{equation} \label{eq:norm:equiv}
   a_{\min}^{1/2} \|\nabla v\|_{L^2(D)} \le
   \| a^{1/2}(\cdot,\y) \nabla v \|_{L^2(D)} \le
   a_{\max}^{1/2} \|\nabla v\|_{L^2(D)}\quad
   \text{$\pi$-a.e. on $\G$}.
\end{equation}
%

The parametric problem~\eqref{eq:pde:strong} is understood in the weak sense:
given $f \in L^2(D)$, find $u : \G \to \X$ such that
\begin{align} \label{eq:pde:weak}
   \int_D a(x, \y) \nabla u(x,\y) \cdot \nabla v(x) \, \dx = \int_D f(x) v(x) \, \dx 
   \quad \forall v \in \X,\ \text{$\pi$-a.e. on $\G$}.
\end{align}
The above assumptions on $a$ and $f$ guarantee that the parametric problem~\eqref{eq:pde:strong}
admits a unique weak solution $u$ in the Bochner space \abrevx{$\V := L_\pi^p(\G; \X)$ for any $p \in [1, \infty]$};
see~\cite[Lemma~1.1]{BabuskaNT07} for details.
In the sequel, we \abrevx{restrict attention to $p=2$ and denote by $\norm{\cdot}{}$ the norm in $\V = L_\pi^2(\G; \X)$;
we also define $\norm{\cdot}{\X} := \norm{\nabla\cdot}{L^2(D)}$}.

The second parametric elliptic problem \eqref{eq:pde2:strong}   combines  \abrevx{uncertainty in the}
source term  with  an isotropic diffusion coefficient field. 
In this case the right-hand side function $f$ simply needs to be a random field that is smooth enough to ensure that
 \eqref{eq:pde2:strong} also admits a unique weak solution $u$
  in the Bochner space $\V$.

\section{Multilevel stochastic collocation finite element method} \label{sec:scfem}

Full details of the construction  of a multilevel stochastic collocation finite element approximation 
of the first parametric elliptic problem can be found  in \partI.  The parametric approximation is associated with
a monotone (or, downward-closed) finite set $\indset_\bullet \subset \N^M$ of multi-indices, where
$\indset_\bullet = \{ \nnu = (\nu_1,\ldots,\nu_M) : \nu_m \in \N, \forall\, m = 1,\ldots,M \}$ 
is such that $\#\indset_\bullet < \infty$. 
\abrevx{Each component $\nu_m$ ($m = 1,\ldots,M$) of the multi-index $\nnu \in \indset_\bullet$ corresponds to}
a set of collocation points along the  $m$th coordinate axis in $\R^M$,  and  the associated {\it sparse grid} 
$\Colpts_\bullet = \Colpts_{\indset_\bullet}$
of collocation points on $\G$ is  given by\footnote{The notation  is 
identical to that  in  \partI. The reader is referred to this paper for any omitted details.}
\[
   \Colpts_{\indset_\bullet} := \bigcup_{\nnu \in \indset_\bullet} \Colpts^{\,(\nnu)}
   = \abrevx{\bigcup_{\nnu \in \indset_\bullet}}\,
      \Colpts_1^{\mf(\nu_1)} \times \Colpts_2^{\mf(\nu_2)} \times \ldots \times \Colpts_M^{\mf(\nu_M)}.
\]
Each collocation point \abrevx{$\z \in \Colpts_{\indset_\bullet} \subset \G$}
is associated with a piecewise linear finite element approximation space 
\abrevx{$\X_{\bullet \z} = \SS^1_0(\TT_{\bullet \z})$} defined on a mesh \abrevx{$\TT_{\bullet \z}$} and
an enhanced space   $\widehat\X_{\bullet \z}$ \abrevx{defined on the mesh $\widehat\TT_{\bullet \z}$}
obtained by {\it uniform refinement}  of \abrev{$\TT_{\bullet \z}$}. 
The spatial detail  space $\Y_{\bullet \z}$ is the approximation space associated with the newly introduced (mid-edge) nodes\abrevx{, i.e.,}
$\widehat \X_{\bullet \z} = \X_{\bullet \z} \oplus \Y_{\bullet \z}$.
\abrevx{We assume that any finite element mesh employed for the spatial discretization is obtained by (uniform or local)
refinement of a given (coarse) initial mesh~$\TT_0$.}

The SC-FEM approximation of the solution $u$ to either of the parametric problems \eqref{eq:pde:strong} or
\eqref{eq:pde2:strong} is given by
\begin{equation} \label{eq:scfem:sol}
   \scsol := \sum\limits_{\z \in \Colpts_\bullet} u_{\bullet \z}(x) \LagrBasis{\bullet \z}{}(\y),
\end{equation}
where $u_{\bullet \z} \in \X_{\bullet \z}$ are Galerkin approximations satisfying~\eqref{eq:sample1:fem} 
or~\eqref{eq:sample2:fem} for $\z \in \Colpts_\bullet$, and
$\{ \LagrBasis{\bullet \z}{}(\y) = \LagrBasis{\z}{\Colpts_\bullet}(\y) : \z \in \Colpts_\bullet \}$ is a set of
multivariable Lagrange basis functions associated with $\Colpts_\bullet$ and satisfying
$\LagrBasis{\bullet \z}{}(\z') = \delta_{\z\z'}$ for any $\z,\,\z' \in \Colpts_\bullet$.
%
The enhancement of the parametric component of the SC-FEM approximation~\eqref{eq:scfem:sol}
is done by enriching the index set~$\indset_\bullet$  
\abrevx{with multi-indices selected} from the {\it reduced margin} set $ \rmarg_{\bullet} = \rmarg({\indset_\bullet})$\abrevx{; this
corresponds to adding \rbl{some} collocation points from the set $\widehat\Colpts_\bullet \setminus \Colpts_\bullet$,
where $\widehat\Colpts_\bullet := \Colpts_{\indset_\bullet \cup \rmarg(\indset_\bullet)}$}.

To keep the discussion concise we simply summarize the components of the adaptive refinement strategy.
The  three components are:
\begin{itemize}
\item 
solution  of a deterministic finite element problem at each sparse grid collocation point. 
That is, the computation of  $u_{\bullet \z} \in \X_{\bullet \z}$ satisfying either
\begin{subequations}
\begin{align} \label{eq:sample1:fem}
   \int_D a(x, \z) \nabla u_{\bullet \z}(x) \cdot \nabla v(x) \, \dx = \int_D f(x) v(x)\, \dx 
   \quad \forall v \in \X_{\bullet \z}
\\
\intertext{in the case of the first parametric problem \eqref{eq:pde:strong}, or} 
 \label{eq:sample2:fem}
   \int_D  \nabla u_{\bullet \z}(x) \cdot \nabla v(x) \, \dx = \int_D f(x, \z) v(x)\, \dx 
   \quad \forall v \in \X_{\bullet \z}
\end{align}
\end{subequations}
in the case of the second parametric problem \eqref{eq:pde2:strong}.
The {\it enhanced} Galerkin solution satisfying~\eqref{eq:sample1:fem} or~\eqref{eq:sample2:fem} 
for all $v \in \widehat\X_{\bullet \z}$
is denoted by $\widehat u_{\bullet \z} \in \widehat\X_{\bullet \z}$.
\item
computation of the spatial hierarchical error \abrevx{indicators}.
For each $\z \in \Colpts_\bullet$, we define $\mu_{\bullet \z} := \norm{e_{\bullet \z}}{\X}$,
where $e_{\bullet \z} \in \Y_{\bullet \z}$ satisfies
\begin{subequations}
\begin{align} \label{eq:hierar1:estimator}
\begin{split}
   \int_D \nabla e_{\bullet \z}(x) \cdot \nabla v(x) \,\dx &=
   \int_D f(x) v(x) \,\dx
   \\
   &\quad -
   \int_D a(x, \z) \nabla u_{\bullet \z}(x) \cdot \nabla v(x) \,\dx \quad \forall v \in \Y_{\bullet \z}
 \end{split}
   \\
   \intertext{in the case of the first parametric problem \eqref{eq:pde:strong}, or satisfies} 
    \label{eq:hierar2:estimator}
  \begin{split}
     \int_D \nabla e_{\bullet \z}(x) \cdot \nabla v(x) \,\dx &=
   \int_D f(x, \z) v(x) \,\dx 
   \\
   &\quad -
   \int_D  \nabla u_{\bullet \z}(x) \cdot \nabla v(x) \,\dx \quad \forall v \in \Y_{\bullet \z}
   \end{split}
\end{align}
\end{subequations}
in the case of the second parametric problem \eqref{eq:pde2:strong}\abrevx{; the corresponding \emph{local}
error indicators $\mu_{\bullet \z}(\xi)$ associated with interior edge midpoints $\xi \in \NN_{\bullet \z}^{+}$
are given by components of the solution vector to the linear system stemming from
the discrete formulation~\eqref{eq:hierar1:estimator} or~\eqref{eq:hierar2:estimator}.}
\item
computation of \abrevx{the 
parametric} error indicators\footnote{This construction assumes that
the enriched index set $\widehat\indset_\bullet$ is obtained using the reduced margin of $\indset_\bullet$, see
Remark~2 in \partI.}
\begin{align} \label{eq:param:indicators:1}
   \widetilde\tau_{\bullet \nnu} =
   \sum\limits_{\z' \in \widetilde\Colpts_{\bullet \nnu}}
   \abrevx{\norm[\bigg]{u_{0 \z'} - \sum\limits_{\z \in \Colpts_\bullet} u_{0 \z} \LagrBasis{\bullet \z}{}(\z')}{\X}} \,
   \norm{\LagrBasisHat{\bullet \z'}{}}{L_\pi^{\abrevx{2}}(\G)}\quad
   \forall \nnu \in \rbl{\rmarg({\indset_\bullet})}\abrevx{,}
\end{align}
where
\abrevx{$\widetilde\Colpts_{\bullet \nnu} \subset \widehat\Colpts_{\bullet} \setminus \Colpts_{\bullet}$
are the collocation points `generated' by the multi-index $\nnu \in  \rmarg({\indset_\bullet})$,
the functions $u_{0 \z'} \in \X_{0 \z'}$ for $\z' \in \widetilde\Colpts_{\bullet \nnu}$
and $u_{0 \z} \in \X_{0 \z}$ for  $\z \in \Colpts_{\bullet}$ are
Galerkin approximations on some meshes $\TT_{0 \z'}$ and  $\TT_{0 \z}$, respectively, that are to be specified
(e.g., $u_{0 \z}$ satisfies~\eqref{eq:sample1:fem} or~\eqref{eq:sample2:fem} with $\X_{\bullet \z}$ replaced by
$\X_{0 \z}$), and}
$\LagrBasisHat{\bullet \z'}{}(\y) = \LagrBasis{\z'}{\widehat\Colpts_\bullet}(\y)$ denotes
the Lagrange polynomial basis function associated with the point $\z' \in \widehat\Colpts_\bullet$  satisfying
$\LagrBasisHat{\bullet \z'}{}(\z'') = \delta_{\z'\z''}$ for any $\z',\,\z'' \in \widehat\Colpts_\bullet$.
\end{itemize}

\abrevx{
We emphasize that the computation of parametric error indicators according to~\eqref{eq:param:indicators:1}
is in line with the hierarchical a posteriori error estimation strategy developed in~\partI\ (see section~4 therein).
In the standard \emph{single-level} SC-FEM setting discussed in~\cite[section~5]{bsx21},
the meshes $\TT_{0 \z'}$ and $\TT_{0 \z}$ underlying the Galerkin approximations $u_{0 \z'}$ and $u_{0 \z}$ in~\eqref{eq:param:indicators:1}
are all selected to be identical to the (single) finite element mesh $\TT_{\bullet \z} = \TT_{\bullet}$ that underlies
the current SC-FEM solution $\scsol$ in~\eqref{eq:scfem:sol}.
In this case, the indicators in~\eqref{eq:param:indicators:1} are written as
\begin{align*} 
   \widetilde\tau_{\bullet \nnu} =
   \sum\limits_{\z' \in \widetilde\Colpts_{\bullet \nnu}}
   \norm{u_{\bullet \z'} - \scsol(\cdot,\z')}{\X} \, \norm{\LagrBasisHat{\bullet \z'}{}}{L_\pi^{\abrevx{2}}(\G)}\quad
   \forall \nnu \in  \rmarg({\indset_\bullet}),
\end{align*}
where $u_{\bullet \z'} \in \X_{\bullet \z'} = \SS^1_0(\TT_\bullet)$ for all $\z' \in \widetilde\Colpts_{\bullet \nnu}$ and
for all $\nnu \in  \rmarg({\indset_\bullet})$.

In the multilevel SC-FEM setting presented in the adaptive algorithm below,
the meshes underlying Galerkin approximations for different collocation points might be different.
In this case, when computing the parametric error indicators in~\eqref{eq:param:indicators:1},
the meshes $\TT_{0 \z'}$~($\z' \in \widetilde\Colpts_{\bullet \nnu}$) and $\TT_{0 \z}$~($\z \in \Colpts_{\bullet}$)
are all selected to be identical to the \emph{coarsest} finite element mesh~$\TT_0$.
}

With \abrevx{the above} ingredients in place, the solution  to the problems in~section~\ref{sec:problem} can be 
generated using  the iterative strategy described in Algorithm~\ref{algorithmx} together with
the marking strategy in Algorithm~\ref{algorithmm}.

\begin{algorithm} \label{algorithmx}
{\bfseries Input:}
$\indset_0 = \{ \1 \}$; 
\abrevx{$\TT_{0 \z} := \TT_0$} for all $\z \in \widehat\Colpts_0 = \Colpts_{\indset_0 \cup \rmarg(\indset_0)}$;
marking criterion.\\
Set the iteration counter $\ell := 0$,  the output counter $k$ and the \rbl{tolerance}.
\begin{itemize}
\item[\rm(i)] 
Compute Galerkin approximations $\big\{ u_{\ell \z} \in \X_{\ell \z} : \z \in \widehat\Colpts_\ell \big\}$ by solving~\eqref{eq:sample1:fem} or~\eqref{eq:sample2:fem}.
\item[\rm(ii)] 
Compute spatial error indicators \abrevx{$\big\{ \mu_{\ell\z} = \norm{e_{\ell \z}}{\X} : \z \in \Colpts_{\ell} \big\}$}
by solving  \eqref{eq:hierar1:estimator}  or \eqref{eq:hierar2:estimator}.
\item[\rm(iii)] 
Compute the parametric error indicators
$\big\{ \widetilde\tau_{\ell \nnu} :  \nnu \in \rbl{\rmarg({\indset_\ell)}}  \big\}$
given by~\eqref{eq:param:indicators:1}.
\item[\rm(iv)] 
Use a marking criterion
to determine $\MM_{\ell \z} \subseteq \NN_{\ell \z}^+$ for all $\z \in \Colpts_\ell$ and
$\markindset_\ell \subseteq  \rbl{\rmarg({\indset_\ell)}}$.
\item[\rm(v)] For all $\z \in \Colpts_\ell$, set $\TT_{(\ell+1) \z} := \refine(\TT_{\ell \z},\MM_{\ell \z})$.
\item[\rm(vi)] Set $\indset_{\ell+1} := \indset_\ell \cup \markindset_\ell$,
%
%
\abrevx{run Algorithm~{\rm \ref{meshalgorithm}} for each
$\z' \in \mathop{\cup}\limits_{\nnu \in \markindset_\ell} \widetilde\Colpts_{\ell \nnu}$
to \rbl{construct meshes} $\TT_{(\ell+1) \z'}$ and
\rbl{initialize} $\TT_{(\ell+1) \z} := \TT_{0 \z} = \TT_0$ for all $\z \in \widehat\Colpts_{\ell+1} \setminus \Colpts_{\ell+1}$}.
\item[\rm(vii)] 
If $\ell = j k$, $j\in \N$, compute the spatial and parametric  {error estimates}  $\mu_\ell$  and  $\tau_\ell$
\abrevx{given by~\eqref{eq:estimate:3} and~\eqref{eq:estimate:8}, respectively,}
and exit if  $\mu_\ell + \tau_\ell < {\tt error tolerance}$.
\item[\rm(viii)] Increase the counter $\ell \mapsto \ell+1$ and goto {\rm(i)}.
\end{itemize}
{\bfseries Output:} For some specific  $\ell_*=jk \in \N$,
the algorithm returns the multilevel SC-FEM approximation~$u_{\ell_*}^{\rm SC}$
computed via~\eqref{eq:scfem:sol} from Galerkin approximations 
$\big\{ u_{{\ell_*} \z} \in \X_{{\ell_*} \z} : \z \in \Colpts_\ell \big\}$
together with a corresponding \rbl{error estimate} $\mu_{\ell_*} + \tau_{\ell_*}$.
\end{algorithm}

A general marking strategy for step~(iv) of Algorithm~\ref{algorithmx} is specified next. We will
adopt this strategy in the numerical experiments discussed in the next section.

\begin{algorithm} \label{algorithmm}
\textbf{Input:}
\abrevx{error indicators}
$\{ \mu_{\ell \z} 
: \z \in \Colpts_\ell \}$,
$\{ \mu_{\ell \z}(\xi) : \z \in \Colpts_\ell,\, \xi \in \NN_{\ell \z}^+ \}$,
and
$\{ \abrevx{\widetilde\tau_{\ell \nnu}} : \nnu \in  \rbl{\rmarg({\indset_\ell)}}  \}$;
marking parameters $0 < \theta_\X, \theta_\Colpts \le 1$ and $\vartheta > 0$.
\begin{itemize}
\item[$\bullet$]
If \
$\sum_{\z \in \Colpts_\ell} \mu_{\ell \z} \norm{L_{\ell \z}}{L^{\abrevx{2}}_{\pi}(\G)}
  \ge \vartheta \sum_{\nnu \in \rbl{\rmarg({\indset_\ell)}} } \abrevx{\widetilde\tau_{\ell \nnu}}$,
then proceed as follows:
\begin{itemize}
\item[$\circ$]
set $\abrev{\markindset_\ell} := \emptyset$
\item[$\circ$]
for each $\z \in \Colpts_\ell$,
determine $\MM_{\ell \z} \subseteq \NN_{\ell \z}^+$ 
\abrevx{such that}
\abrevx{
\begin{equation} \label{eq:doerfler:separate1}
 \theta_\X \, \sum_{\z \in \Colpts_\ell} \sum_{\xi \in \NN_{\ell \z}^+} \mu_{\ell \z}(\xi) \norm{L_{\ell \z}}{L^{\abrevx{2}}_{\pi}(\G)} \le
 \sum_{\z \in \Colpts_\ell} \sum_{\xi \in \MM_{\ell \z}} \mu_{\ell \z}(\xi) \norm{L_{\ell \z}}{L^{\abrevx{2}}_{\pi}(\G)}
\end{equation}
\rbl{with a cumulative cardinality $\sum_{\z \in \Colpts_\ell} \#\MM_{\ell \z}$ that is  minimized
over all the sets that satisfy~\eqref{eq:doerfler:separate1}}.
}
\end{itemize}
\item[$\bullet$]
Otherwise, if \ 
$\sum_{\z \in \Colpts_\ell} \mu_{\ell \z} \norm{L_{\ell \z}}{L^{\abrevx{2}}_{\pi}(\G)}
  < \vartheta \sum_{\nnu \in \rbl{\rmarg({\indset_\ell)}} } \abrevx{\widetilde\tau_{\ell \nnu}}$,
then proceed as follows:
\begin{itemize}
\item[$\circ$]
set $\MM_{\ell \z} := \emptyset$ for all $\z \in \Colpts_\ell$
\item[$\circ$]
determine $\markindset_\ell \subseteq \rbl{\rmarg({\indset_\ell)}} $
of minimal cardinality such that
\begin{equation} \label{eq:doerfler:separate2}
 \theta_\Colpts \, \sum_{\nnu \in \rbl{\rmarg({\indset_\ell)}} } \abrevx{\widetilde\tau_{\ell \nnu}} \le
 \sum_{\abrev{\nnu \in \markindset_\ell}} \abrevx{\widetilde\tau_{\ell \nnu}}.
\end{equation}
\end{itemize}
\end{itemize}
\textbf{Output:}
$\MM_{\ell \z} \subseteq \NN_{\ell \z}^+$ for all $\z \in \Colpts_\ell$ and
$\markindset_\ell \subseteq  \rbl{\rmarg({\indset_\ell)}}$.
\end{algorithm}

As discussed in section~4 of \partI, the
computation of the {\it error estimates}  in step~(vii) of Algorithm~\ref{algorithmx} is best done
periodically because of the significant computational overhead. 
Specifically, the spatial error estimate  
\begin{align}  \label{eq:estimate:3}
   \mu_\bullet   & := 
    \norm[\bigg]{\sum\limits_{\z \in \Colpts_\bullet} (\widehat u_{\bullet \z} - u_{\bullet \z})\, \LagrBasis{\bullet \z}{}}{}
\end{align}
requires computation of the enhanced Galerkin approximation \abrevx{$\widehat u_{\bullet \z} \in \widehat\X_{\bullet \z}$}
and thus requires the solution of the PDE on
\abrevx{the mesh $\widehat\TT_{\bullet \z}$---a uniform refinement of $\TT_{\bullet \z}$---for each collocation point}
generated by the current index set.
The parametric error estimate \abrev{(cf.~\eqref{eq:param:indicators:1})} 
%
%
%
%
\begin{align}
 \label{eq:estimate:8}
   \tau_\bullet  & : =
   \norm[\bigg]
   {\sum\limits_{\z' \in \widehat\Colpts_\bullet \setminus \Colpts_\bullet}
                        \abrevx{ \Big( u_{0 \z'} - \sum\limits_{\z \in \Colpts_\bullet} u_{0 \z} \LagrBasis{\bullet \z}{}(\z') \Big) }
                        \LagrBasisHat{\bullet \z'}{}}{}
\end{align}
requires\abrevx{, as discussed above,} additional PDE solves \abrev{on the coarsest mesh~$\TT_{0 \z'} := \TT_0$ for all}
margin collocation points $\z' \in \widehat\Colpts_\bullet \setminus \Colpts_\bullet$
\abrevx{(the coarsest-mesh Galerkin approximations $u_{0 \z}$ in~\eqref{eq:estimate:8}
for the current collocation points $\z \in \Colpts_\bullet$ will have been computed in preceding iterations 
and, thus, can be \rbl{reused})}.
The key point here  is that computation of the \abrevx{error} estimates is only needed
to give  reliable termination of the adaptive process (and  to provide reassurance 
that the SC-FEM error is decreasing at an acceptable rate).

\abrevx{Regarding the implementation aspects of computing the above error estimates, we
note that the sum in~\eqref{eq:estimate:3} involves Galerkin approximations over different finite element meshes.
In our implementation, the computation of this sum is effected by
interpolating piecewise linear functions $u_{\bullet \z}$ and $\widehat u_{\bullet \z}$
at the nodes of the mesh 
$\bigoplus_{\z \in \Colpts_\bullet} \widehat\TT_{\bullet \z}$---the overlay (or, the coarsest common refinement) of the meshes
$\widehat\TT_{\bullet \z}$, $\z \in \Colpts_\bullet$---and by subtracting/summing the obtained coefficient vectors
representing these piecewise linear functions over the same mesh
$\bigoplus_{\z \in \Colpts_\bullet} \widehat\TT_{\bullet \z}$.
In this respect, the implementation of the parametric error estimate in~\eqref{eq:estimate:8} is rather straightforward,
as the involved Galerkin approximations $u_{0 \z}$ and $u_{0 \z'}$ are all computed on the same coarsest finite element mesh~$\TT_0$.
}

The other detail that is missing in the statement of Algorithm~\ref{algorithmx} 
is the identification of a strategy  for defining suitable meshes \abrevx{$\TT_{(\ell+1), \z'}$}
corresponding to the \abrevx{newly `activated'} collocation  points
in step~(vi).  
This specification of sample-specific {\it initial meshes}  turns out to be
crucial if \rbl{optimal rates of convergence are to be realized in practice}.
If an initial mesh \abrevx{associated with a collocation point} is too coarse, then
\abrevx{`activating' this} collocation point will introduce a large spatial error at the next iteration step.
Conversely, if the initial mesh is too fine, as in the case
of a single-level implementation of the algorithm, then the  growth in the \abrevx{number of} 
degrees of freedom is not matched by the resulting error reduction. Indeed, the conclusion 
reached in~\cite{FeischlS21} on this point
is that ``while the theoretical results are strongest for the fully adaptive algorithm ... the 
single mesh algorithm seems to be more efficient''.
A  mesh initialization  strategy that attempts to  balance the  conflicting requirements is 
given in Algorithm~\ref{meshalgorithm}.
\abrevx{Specifically, for a given (newly `activated') collocation point $\z' \not\in \Colpts_{\bullet}$,
we start with the coarsest mesh $\TT_0$
and iterate the standard SOLVE $\to$ ESTIMATE $\to$ MARK $\to$ REFINE loop until the resolution of the mesh is such that
the estimated error in the corresponding Galerkin solution $u_{\bullet \z'}$ is on par with the error estimates for Galerkin solutions
associated with other (already `active') collocation points $\z \in \Colpts_{\bullet}$.
This is ensured by the choice of stopping tolerance ${\tt tol}$ in Algorithm~\ref{meshalgorithm}.
We note that in the multilevel SGFEM, such a mesh initialization procedure is not needed.
Instead, for every newly `activated' multi-index,
the associated finite element mesh is set to the coarsest mesh $\TT_0$; see~\cite{bpr2020+}.
Due to the inherent orthogonality of the parametric components of SGFEM approximations associated with different multi-indices,
this initialization by the coarsest mesh does not affect optimal convergence properties of the multilevel SGFEM; see~\cite{bpr2021+}.
}

\begin{algorithm} \label{meshalgorithm}
\textbf{Input:}
spatial error indicators \abrevx{$\big\{ \mu_{\ell\z}: \z \in \Colpts_{\ell} \big\}$};
\rblx{the set of collocation points $\Colpts_{\ell+1} = \Colpts_{\indset_{\ell+1}}$;}
the collocation point~\rblx{$\z' \in \Colpts_{\ell+1} \setminus \Colpts_{\ell}$};
marking parameter $\theta$.\\
\abrevx{
Set the tolerance
${\tt tol} := (\# \Colpts_{\ell})^{-1} \sum_{\z \in \Colpts_{\ell}} \mu_{\ell \z} \norm{L_{\rblx{(\ell+1)} \z}}{L^{2}_{\pi}(\G)}$
and the iteration counter $n := 0$;
initialize the mesh $\TT_{0 \z'} := \TT_0$.
\begin{itemize}
\item[\rm(i)] 
Compute the Galerkin approximation $u_{n \z'} \in \X_{n \z'}$ by solving~\eqref{eq:sample1:fem} or~\eqref{eq:sample2:fem}.
\item[\rm(ii)] 
Compute the error estimate $\mu_{n \z'} = \norm{e_{n \z'}}{\X}$
by solving  \eqref{eq:hierar1:estimator}  or \eqref{eq:hierar2:estimator} and compute the corresponding
local error indicators $\big\{ \mu_{n \z'}(\xi): \xi \in \NN^{+}_{n \z'} \big\}$.
\item[\rm(iii)] 
If $\mu_{n \z'} \norm{L_{\rblx{(\ell+1)} \z'}}{L^{2}_{\pi}(\G)} < {\tt tol}$, set $\TT_{(\ell+1) \z'} := \TT_{n \z'}$ and exit.
\item[\rm(iv)]
Determine $\MM_{n \z'} \subseteq \NN_{n \z'}^+$ of minimal cardinality such that
\begin{equation*} 
   \theta \, \sum_{\xi \in \NN_{n \z'}^+} \mu_{n \z'}(\xi)^2 \le
   \sum_{\xi \in \MM_{n \z'}} \mu_{n \z'}(\xi)^2.
\end{equation*}
\item[\rm(v)] Set $\TT_{(n+1) \z'} := \refine(\TT_{n \z'},\MM_{n \z'})$.
\item[\rm(vi)] Increase the counter $n \mapsto n+1$ and goto {\rm(i)}.
\end{itemize}
}
\textbf{Output:}
\rbl{The mesh} $\TT_{(\ell+1) \z'}$ 
\abrevx{associated with} the collocation point $\z'$.
\end{algorithm}

Results presented in the next section will
show that a well-designed multilevel strategy can give significant efficiency gains 
compared to a \rblx{single-level SC-FEM} algorithm if the parameterized problem has local features 
that \abrevx{vary} in spatial location across the parameter~space.

\section{Numerical experiments}\label{sec:results}

Results for  three test cases are discussed in this section of the paper. 
The performance of our adaptive SC multilevel algorithm will be directly compared with that of the 
single-level algorithm discussed in  \partI\ to see if any gains in efficiency 
can be realized.
The first two test cases are identical to those discussed in \S5 of \partI. The third test case is  a refinement of 
the {\it one peak} test problem that was introduced by Kornhuber \& Youett~\cite{ky18} in order to assess
the efficiency of adaptive Monte Carlo methods.

The single-level refinement strategy that is  the basis for comparison  is the obvious and
natural   simplification of the multilevel strategy described in \S\ref{sec:scfem}. 
Thus, at each step $\ell$ of the process, 
we compute the error indicators associated with the SC-FEM solution $u_{\ell\z}$
(steps (ii)--(iii) of Algorithm~\ref{algorithmx}).
 The marking criterion in Algorithm~\ref{algorithmm} then 
 identifies the refinement type by comparing the (global) spatial error estimate
$\bar\mu_{\ell} := \| \mu_{\ell \z} \norm{L_{\ell \z}}{L^2_{\pi}(\G)} \|_{\ell_1}$ with
the  parametric error estimate $\bar\tau_{\ell} := \| \abrevx{\widetilde\tau_{\ell \nnu}} \|_{\ell_1}$.
To effect  a spatial refinement in the single-level case,
we use \abrevx{a D{\" o}rfler-type} marking with \rbl{threshold} $\theta_\X$ to produce sets of marked elements 
from the (single) grid $\TT_\ell$. 
A refined triangulation $\TT_{\ell+1}$ can then be constructed by refining the elements in the {\it union} 
of these individual sets $\MM_{\ell \z}$ ($\z \in \Colpts_\ell$) of marked elements.

\subsection{Test case I:  affine coefficient data}\label{sec:affineresults}

We set $f = 1$ and look to solve the first model problem on the square-shaped 
domain $D = (0, 1)^2$ with random field coefficient given by
\begin{align} \label{kl}
a(x, \y) = a_0(x) + \sum_{m = 1}^M a_m(x) \, y_m,\quad
   x \in D,\ \y \in \Gamma.
\end{align}
The  specific problem  we consider is taken  from~\cite{bs16}.
The parameters $y_m$ in \eqref{kl} are the images of uniformly 
distributed independent mean-zero random variables, so that $\pi_m = \pi_m(y_m)$ is
the associated probability measure on $\G_m = [-1,1]$.
The expansion coefficients $a_m$, $m\, \in \, \N_0$  are chosen
to represent planar Fourier modes of increasing total order.
Thus, we fix  $a_0(x) := 1$ and set
\begin{equation}
\label{diff_coeff_Fourier}
a_m(x) := \alpha_m \cos(2\pi\beta_1(m)\,x_1) \cos(2\pi\beta_2(m)\,x_2),\  x=(x_1,x_2) 
\in (0,1) \times (0,1).
\end{equation}
The modes are ordered so that for any $m \in \N$,
\begin{equation}
  \beta_1(m) = m - k(m)(k(m)+1)/2\ \ \hbox{and}\ \ \beta_2(m) =k(m) - \beta_1(m)
\end{equation}
with $k(m) = \lfloor -1/2 + \sqrt{1/4+2m}\rfloor$ and the amplitude coefficients are 
constructed so that
$\alpha_m = \bar\alpha m^{-2}$  with $ \bar\alpha = 0.547$.
This is referred to as the {\it slow decay case} in~\cite{bs16}.

\begin{figure}[!pth]
\centering
\includegraphics[width = 0.45\textwidth]{{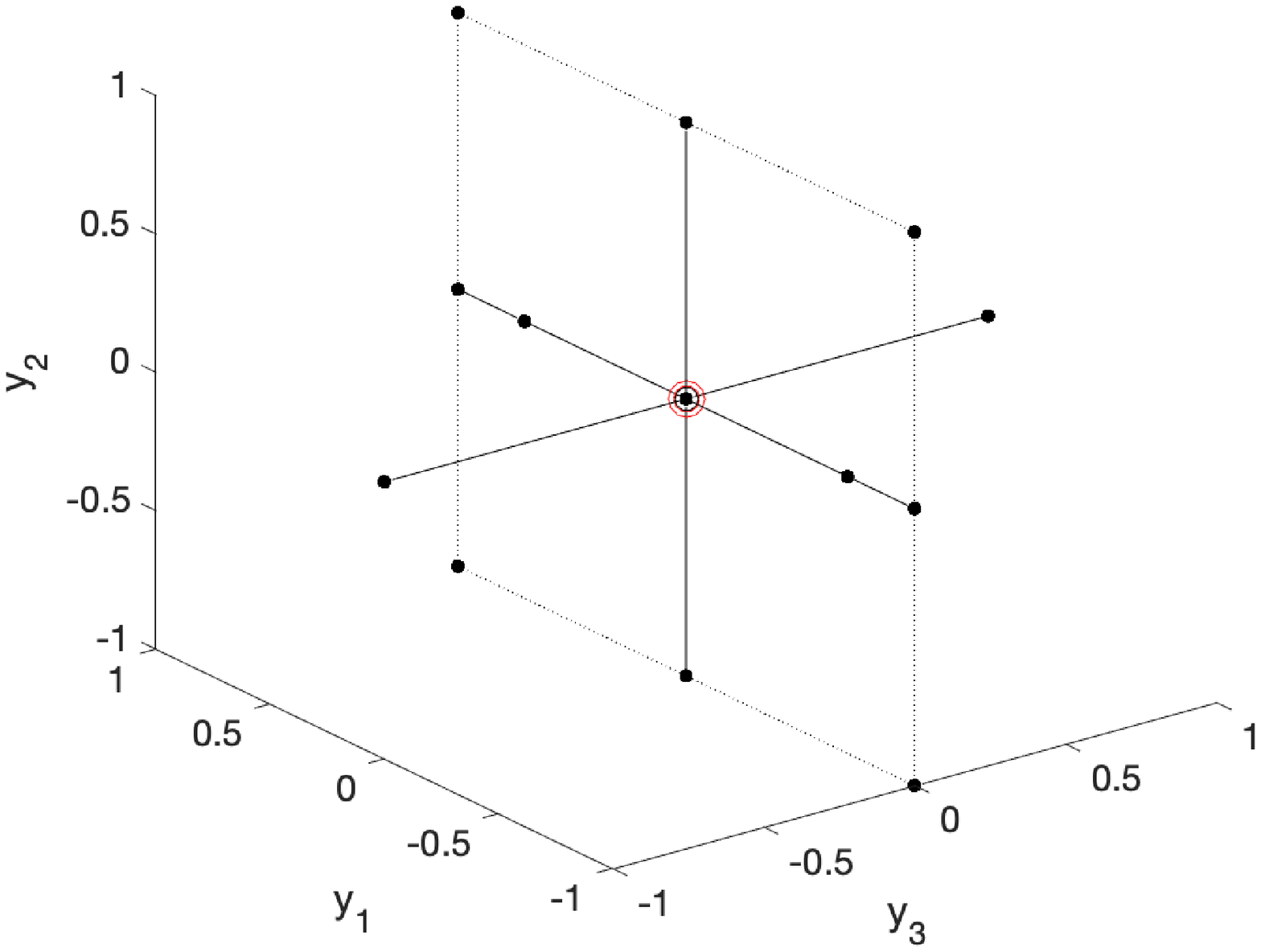}}
\includegraphics[width = 0.45\textwidth]{{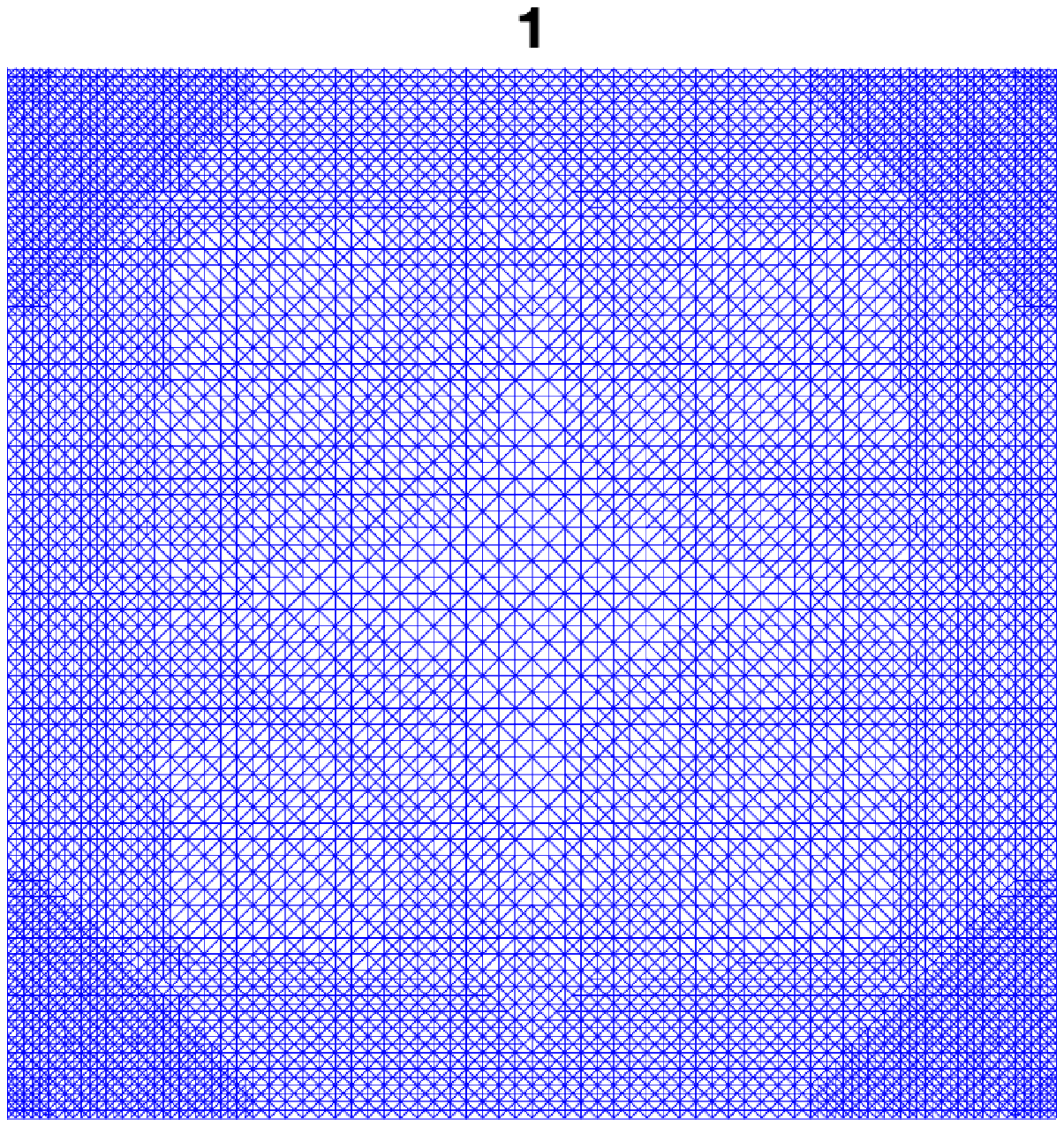}}
\includegraphics[width = 0.45\textwidth]{{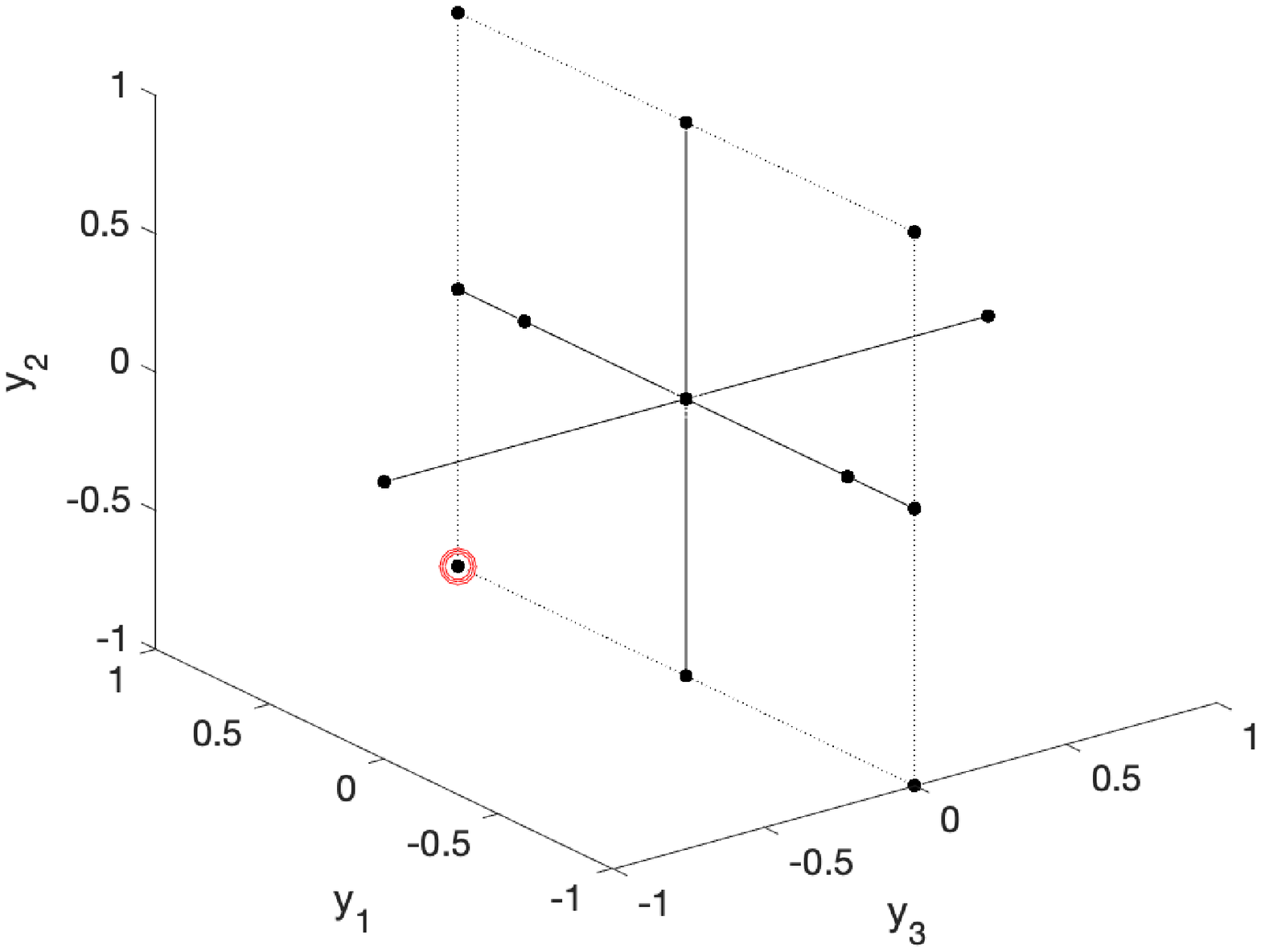}}
\includegraphics[width = 0.45\textwidth]{{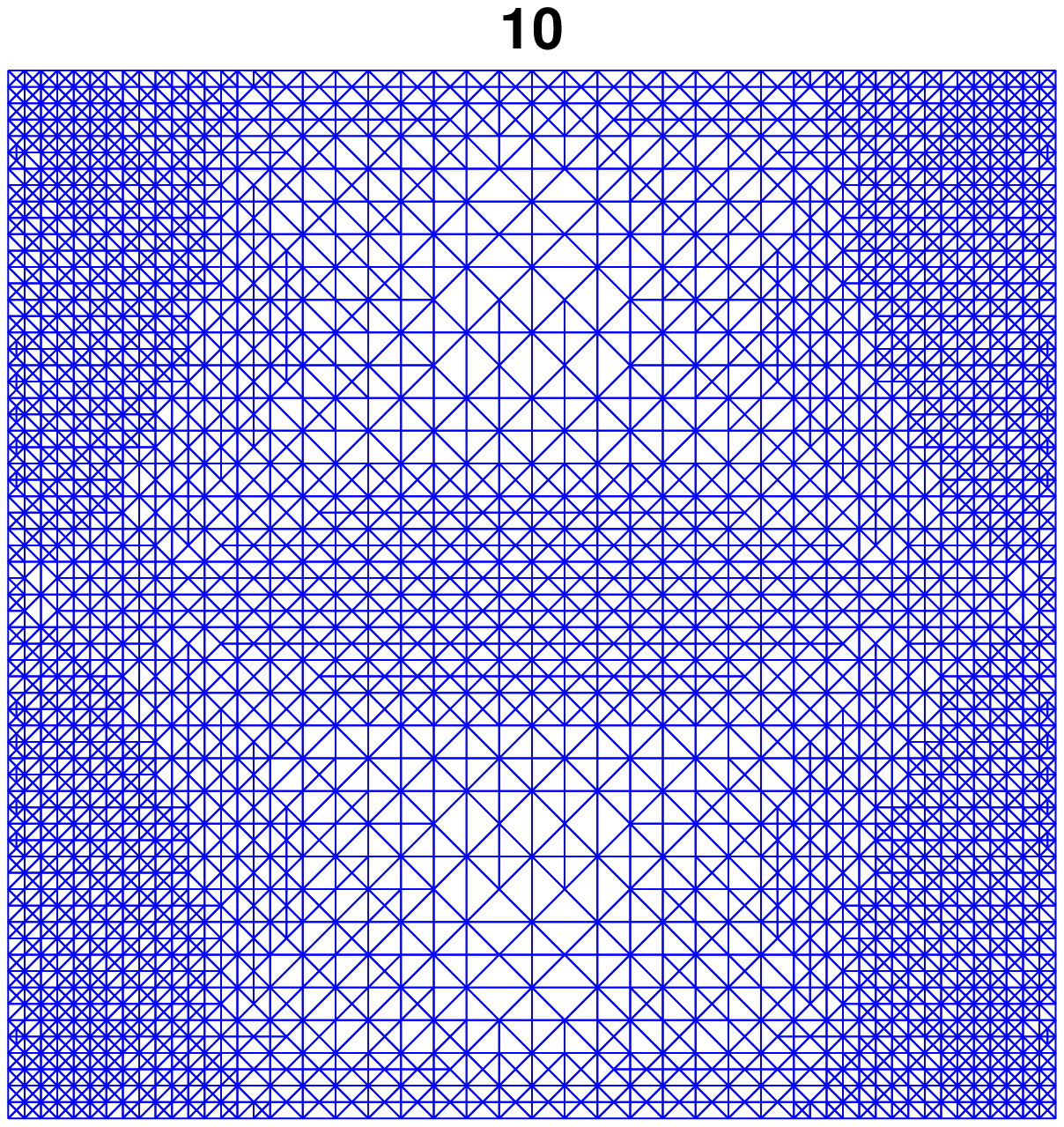}}
\includegraphics[width = 0.45\textwidth]{{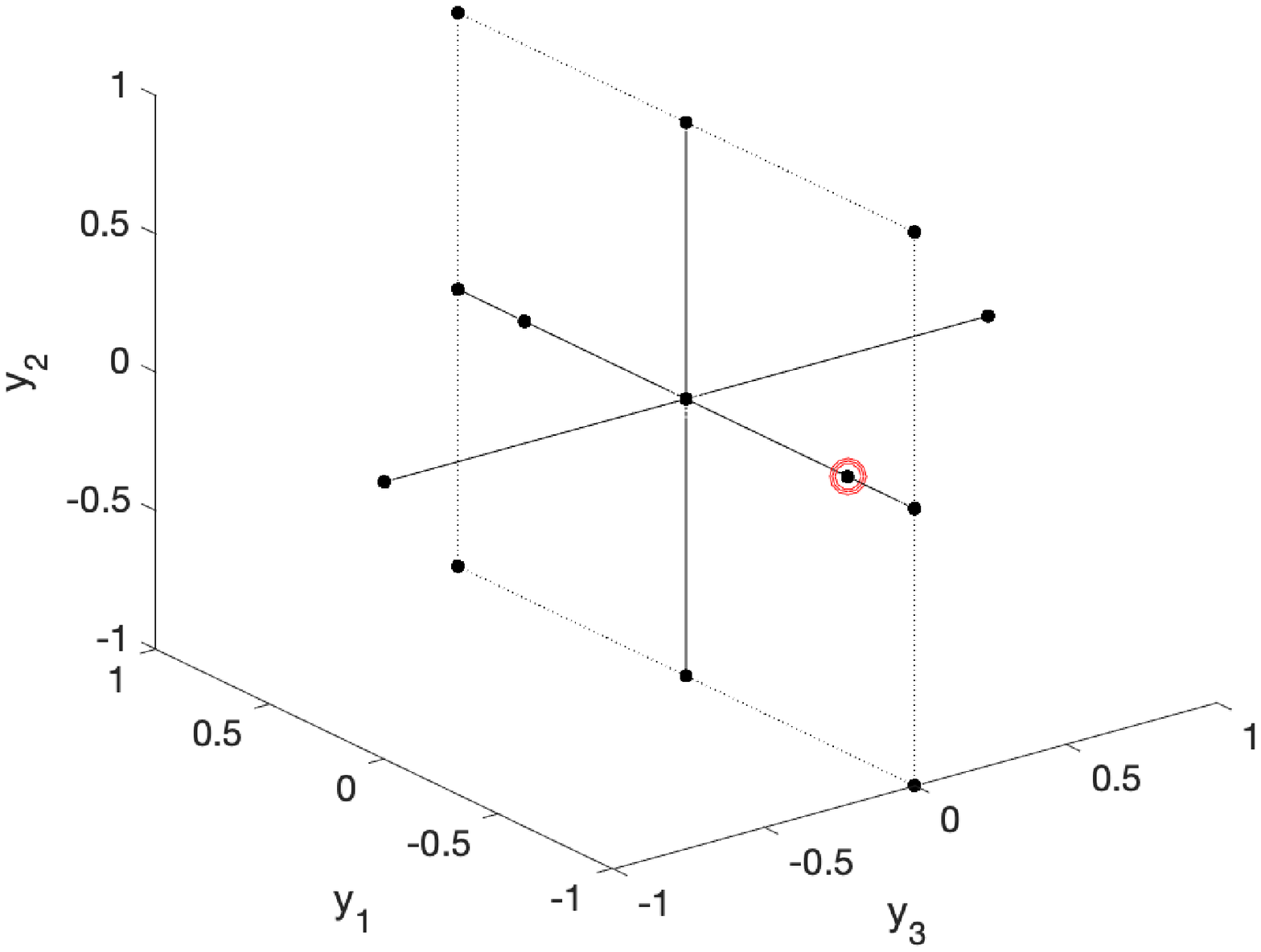}}
\includegraphics[width = 0.45\textwidth]{{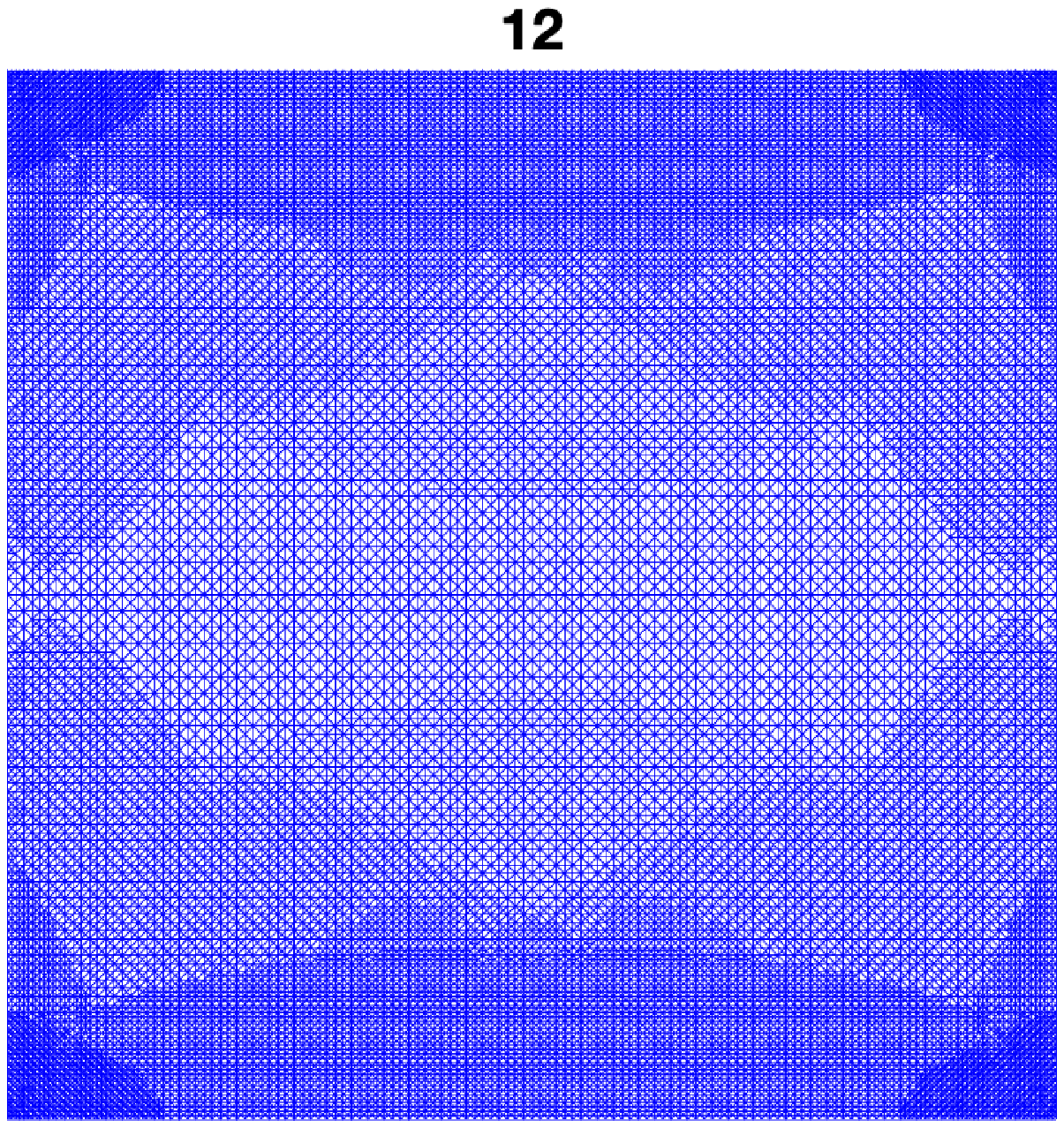}}
\caption{Selected collocation point  (left) and  corresponding spatial mesh (right) that is generated by the 
multilevel adaptive strategy for test case I.}
\label{fig:sc4.1meshes}
\end{figure}

A reference solution to this problem  with $M$ set to 4  is illustrated in Fig.~\rblx{1} in~\partI. 
This  solution  was generated by running the {\it single-level} algorithm with   
the $\tt{error tolerance}$ set to {\tt 6e-3},  starting
from a uniform initial mesh with 81 vertices and a sparse grid 
consisting of a single collocation point. The threshold parameter  $\vartheta$
was set to {\tt 1}, the marking parameters $\theta_\X$ and $\theta_\Colpts$ were set to {\tt 0.3}.
The error tolerance  was satisfied after 25 iterations comprising
 {20} spatial refinement steps and  {5} parametric refinement steps.
 There were {13} Clenshaw--Curtis sparse grid collocation points when 
 the iteration terminated. These  points are visualized in Fig.~\ref{fig:sc4.1meshes}.
 The associated sparse grid indices are  listed in Table~1 in \partI.
The  final spatial mesh  is shown in  Fig.~2 in~\partI. 
The number of vertices in this mesh is {\tt 16,473} so the total number of 
degrees of  freedom when the error tolerance was satisfied when running the 
single-level algorithm was {\tt 214,149}.

\begin{figure}[!thp]
\centering
\includegraphics[width = 0.65\textwidth]{{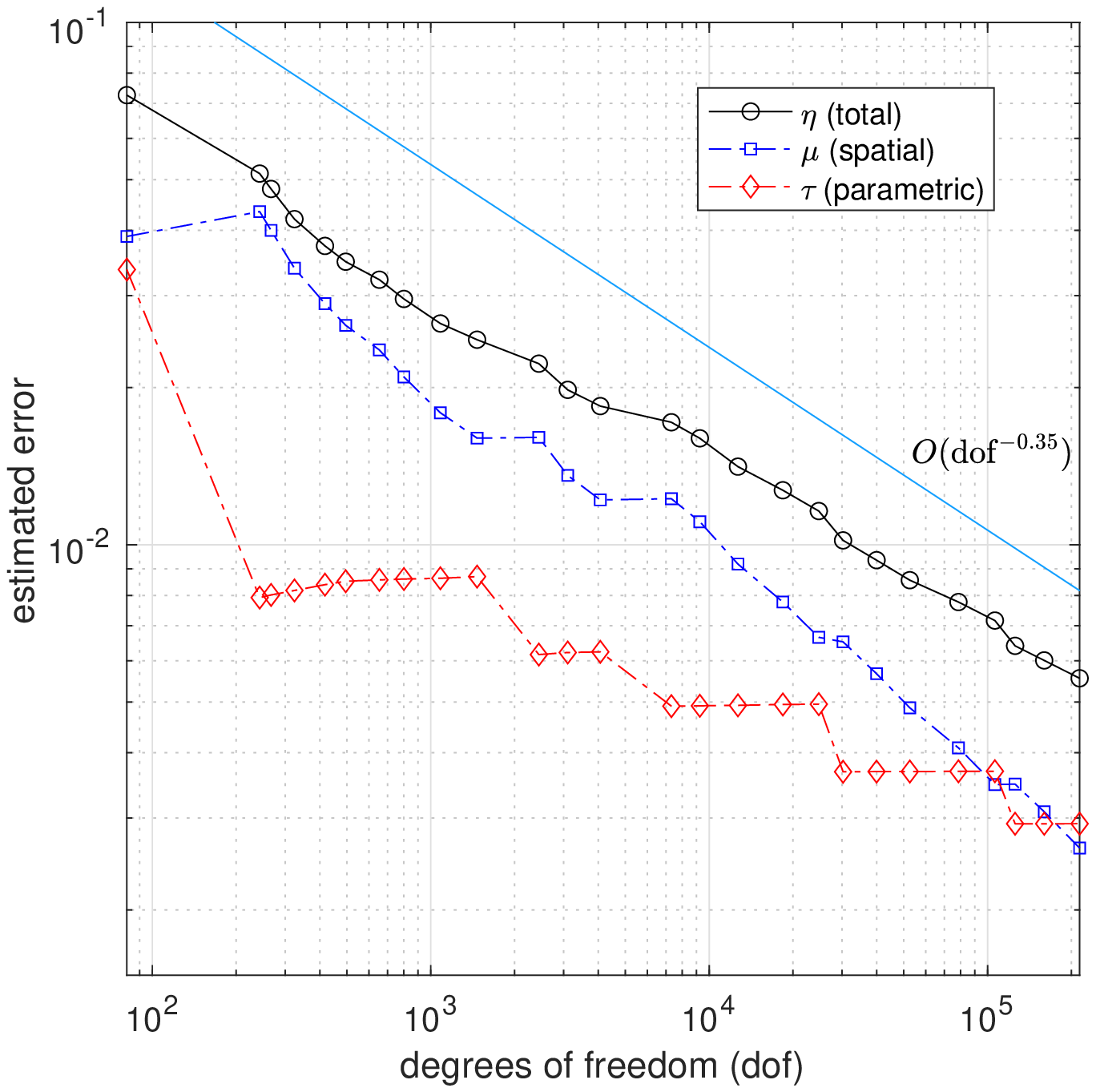}}
\includegraphics[width = 0.65\textwidth]{{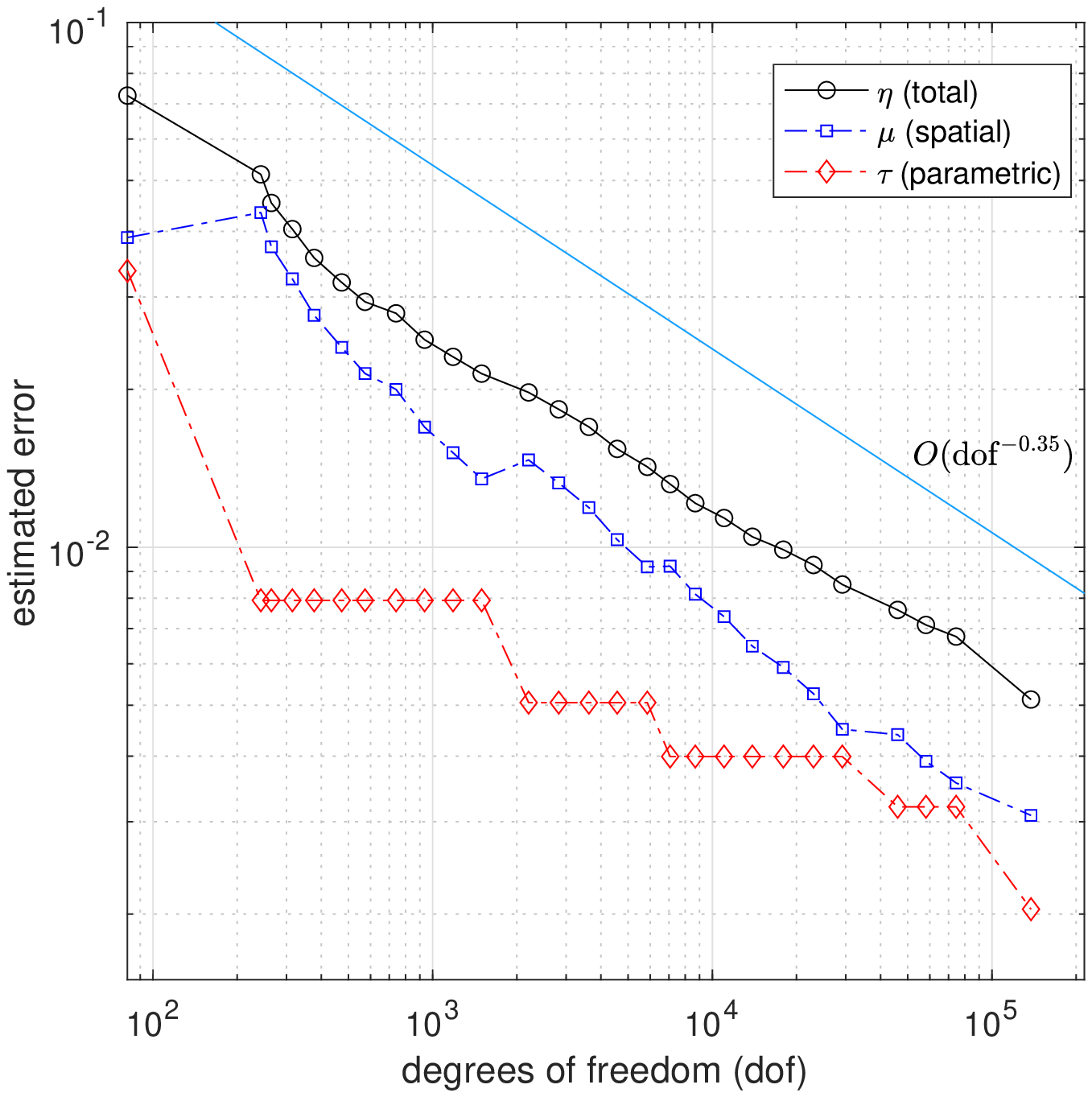}}
\caption{Evolution of the single-level error estimates  (top) and  the
multilevel error estimates  (bottom) 
for test case I with error tolerance set to {\tt 6e-3}.}
\label{fig:sc4.1errors}
\end{figure}

The first test of the {\it multilevel}  algorithm is  to repeat the above experiment;
that is, starting from the same point with identical marking parameters  
$\vartheta=1$, $\theta_\X = \theta_\Colpts = 0.3$
\rblx{(we also set the marking parameter $\theta$ in Algorithm~\ref{meshalgorithm}
to the same value as $\theta_\X$ in all our experiments)}.
 Specifying the same error tolerance {\tt 6e-3} led to the
 the same 13 collocation points being activated, in this case after 26 rather 
 than 25 iterations.
A comparison of the single-level and multilevel error \rblx{estimates} 
is given in Fig.~\ref{fig:sc4.1errors}.  While the final number of degrees of freedom 
is reduced  from  {\tt 214,149} to   {\tt 137,943} in the multilevel case, the 
{\it rate of convergence} is  still far from optimal (close to $O({\rm dof}^{-1/3})$).

The degree of refinement of the final meshes associated with some specific
collocation points is illustrated in Fig.~\ref{fig:sc4.1meshes}. The two finest  meshes
had over {\tt 32,000} vertices and are associated with the pair of collocation points 
that are activated by the sparse grid index {\tt 3 1 1 1} that is introduced  at the 
final  iteration (one of \rblx{these collocation points and the corresponding mesh are} shown in  the bottom plot).
The two coarsest meshes had close to {\tt 3,600} vertices; one of \rblx{these} is shown in the middle plot. 
The \abrevx{mesh} that is associated with the mean field $a_0=1$  has {\tt 11,157}  vertices
and is shown in the topmost plot. As might be anticipated, the level  of refinement 
of this \abrevx{mesh} is less than that of the final \abrevx{mesh} that is generated by the 
single-level~strategy.

\abrevx{It is worth pointing out that in our extensive experimentations with other choices of marking parameters
the adaptive multilevel SC-FEM algorithm did not exhibit a faster convergence rate compared to that of the single-level algorithm
for the respective choice of marking parameters.
This is in contrast to SGFEM, where multilevel adaptivity always results in a faster convergence rate than that of the single-level counterpart
for problems with affine-parametric coefficients including the test case considered here; see~\cite{egsz14, cpb18+, bpr2020+, bpr2021+}.
Furthermore, for this class of problems, the analysis in~\cite{bpr2021+} has shown that,
under an appropriate saturation assumption,
the adaptive multilevel SGFEM algorithm driven by a two-level a posteriori error estimator and employing a D{\" o}rfler-type marking
on the joint set of spatial and parametric indicators yields optimal convergence rates with respect to the number of degrees of freedom
in the underlying multilevel approximation space.
}

\subsection{Test case II:  nonaffine coefficient data}\label{sec:nonaffineresults}

\begin{figure}[!thp]
\centering
\includegraphics[width = 0.65\textwidth]{{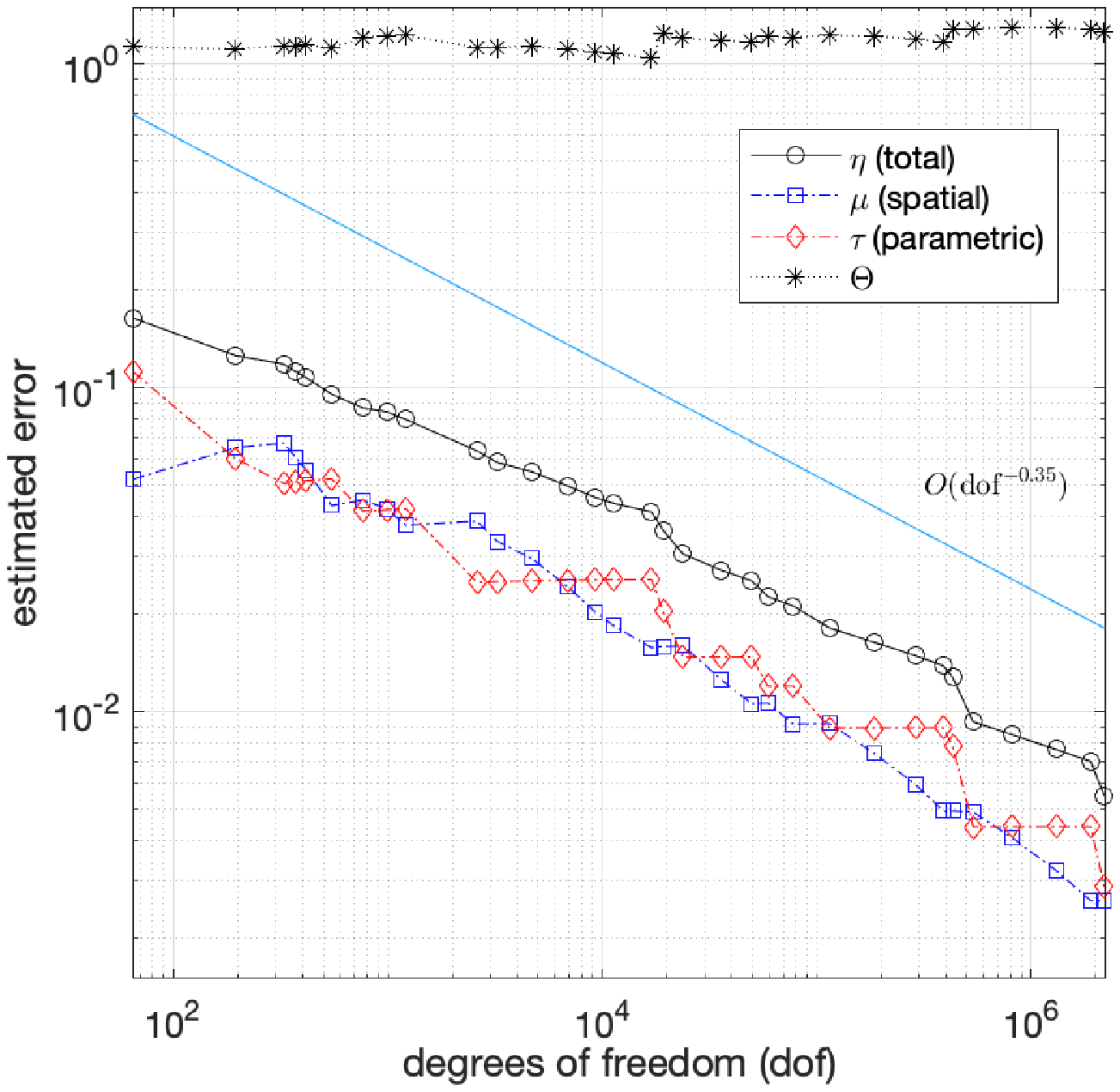}}
\includegraphics[width = 0.65\textwidth]{{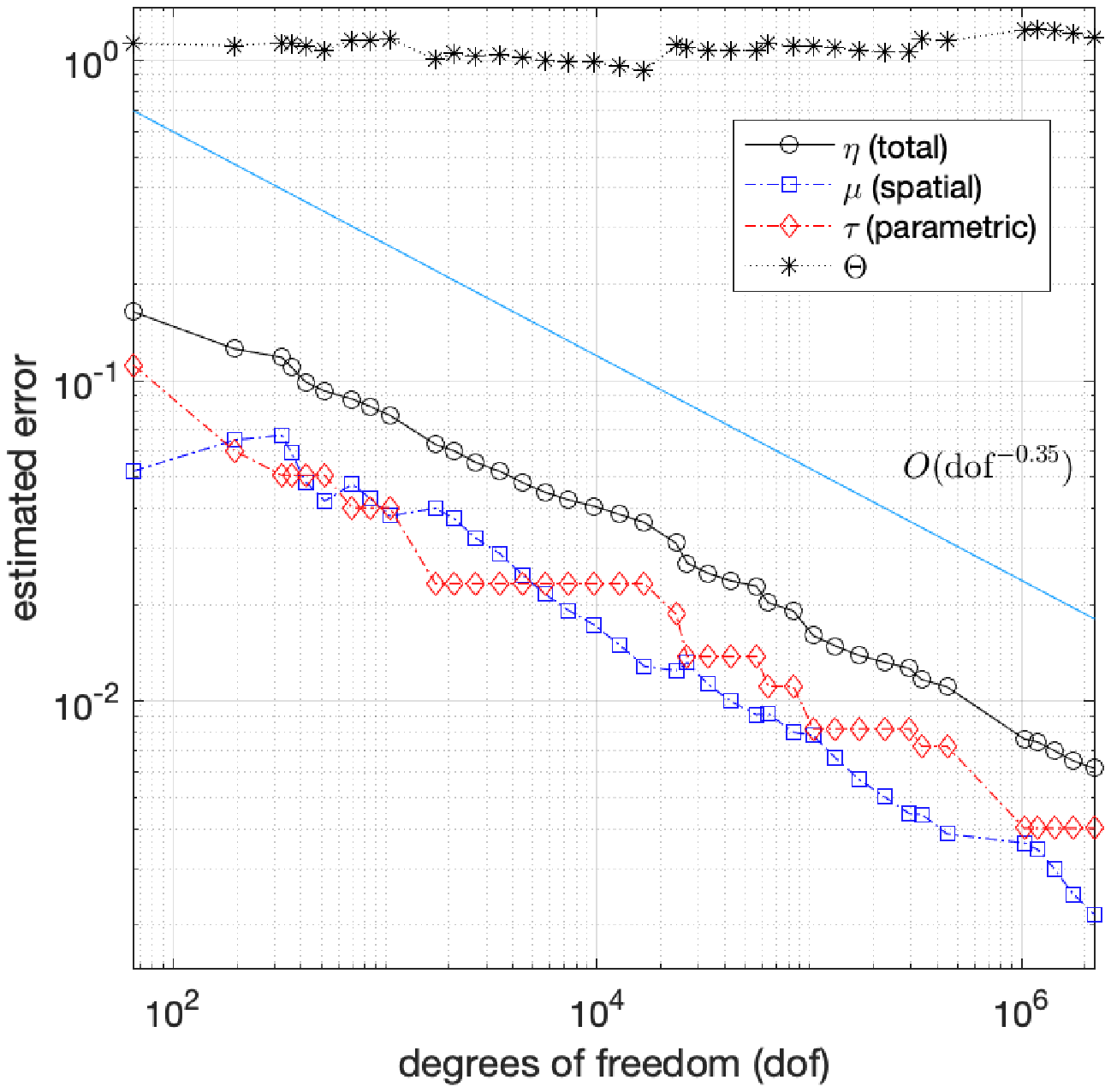}}
\caption{Evolution of the single-level error estimates  (top) and  the
multilevel error estimates  (bottom) 
for test case II with error tolerance set to {\tt 6e-3}.}
\label{fig:sc4.2errors}
\end{figure}

In this case, we set $f = 1$ and look to solve  the first model problem
on the L-shaped domain $D = (-1, 1)^2\backslash (-1, 0]^2$ with
coefficient $a(x, \y) = \exp(h( x, \y))$,
where the exponent field  $h(x, \y)$ has affine dependence on
parameters $y_m$  that are images of uniformly 
distributed independent mean-zero random variables,
\begin{align} \label{kll}
h(x, \y) = h_0(x) + \sum_{m = 1}^{\rbl{4}} h_m(x) \, y_m,\quad
   x \in D,\ \y \in \Gamma.
\end{align}
We further specify $h_0(x) \,{=}\, 1$ and 
$h_m(x) = \sqrt{\lambda_m} \varphi_m(x)$
($m = 1,\ldots, \rbl{4}$). Here $\{(\lambda_m, \varphi_m)\}_{m=1}^\infty$ 
are the eigenpairs of the integral operator
$\int_{\abrevx{ D \cup (-1,0]^2} } \hbox{\rm Cov}[\rblx{h}](x, x') \varphi(x')\, \hbox{d} x' $
with a  synthetic covariance function given by
\begin{align} \label{cov}
\hbox{\rm Cov}[\rblx{h}](x, x') 
= \sigma^2 \exp
\left( -| x_1 - x_1' |  - | x_2 - x_2' |  \right).
\end{align}
The standard deviation  $\sigma$  is set to 1.5 in order
to  mirror   the most challenging test case  in \S5.2 of \partI.
The convergence of the multilevel adaptive algorithm, starting  with one collocation point and
with the initial grid shown in Fig.~7 of \partI\  is compared with the single-level result
in Fig.~\ref{fig:sc4.2errors}.  The multilevel algorithm is again run 
using the marking parameters $\theta_\X = \theta_\Colpts = 0.3$ 
specified in \partI\ and the same error tolerance, that is  {\tt 6e-3}.

These results reinforce the view that performance gains from the multilevel
strategy are  difficult to realize. While the number of active collocation points
is smaller in the multilevel case (51 vs 57; the sparse grid index  {\tt 2  1  2  2}
added at the final single-level iteration is not included), the total number
of degrees of freedom when the tolerance is reached is almost identical
({\tt 2,212,393} vs {\tt 2,190,847}). The issue here is that meshes associated with
mixed indices with multiple active dimensions have multiple  features that
require resolution. Thus, the most refined grid associated with the index that is 
introduced in the final parametric enhancement  has {\tt 428,972} vertices.
This is significantly more refined than the final grid that is generated
in the single-level implementation, which had {\tt 37,133} vertices. This fact,
together with the  increase in the number of adaptive steps taken (37 vs 31) means 
that the overall computation time is significantly increased when the
multilevel strategy is adopted.

\abrevx{The plots in Fig.~\ref{fig:sc4.2errors} also show that the use of the \emph{coarsest-mesh} approximations
for computing the parametric error estimates $\tau_{\ell}$ in~\eqref{eq:estimate:8} does not affect
the overall effectivity of the error estimation in the multilevel algorithm.
Indeed, in the single-level algorithm (where parametric error estimates employ the (single) \emph{refined mesh}
underlying the current SC-FEM solution~$u_{\ell}^{\rm SC}$),
the effectivity indices $\Theta_\ell$ computed\footnote{The effectivity indices are computed using a reference solution
as explained in~\cite{bsx21}, see equation~(42)~therein.}
at each iteration range between 1.047 and 1.296,
whereas for the multilevel algorithm they stay between 0.930 and~1.257. 
}

\subsection{Test case III:  one peak problem}\label{sec:onepeakresults}


We are looking to solve the Poisson equation $-\nabla^2 u = f $ in a unit square domain 
$D=(-4,4)\times(-4,4)$ with Dirichlet boundary data  $u=g$.   
The source term $f$ and boundary data are  {\it uncertain} and  are
parameterized by  $ \y =(y_1,y_2)$, representing the image  of  a pair 
of independent random  variables  with  $y_j \sim {U}[-1,1]$. 
In the vanilla case discussed in~\cite{ky18}, the same test problem is \rblx{posed} 
on  the unit domain $I=(-1,1)\times(-1,1)$ with  $y_j \sim {U}[-1/4,1/4]$. 
The source term $f$ and the  boundary data $g$ are chosen so that the  problem 
has a specific pathwise solution given by 
\begin{align} \label{peaksolv}
u(x, \y) & = \exp ( - \beta  \{ (x_1 -y_1)^2 + (x_2 -y_2)^2 \} ), 
\end{align} 
where a scaling factor $\beta=50$ is chosen to generate
a highly localized Gaussian profile centered at the  uncertain spatial location $(y_1,y_2)$.

\begin{figure}[!thp]
\centering
\includegraphics[width = 0.8\textwidth]{{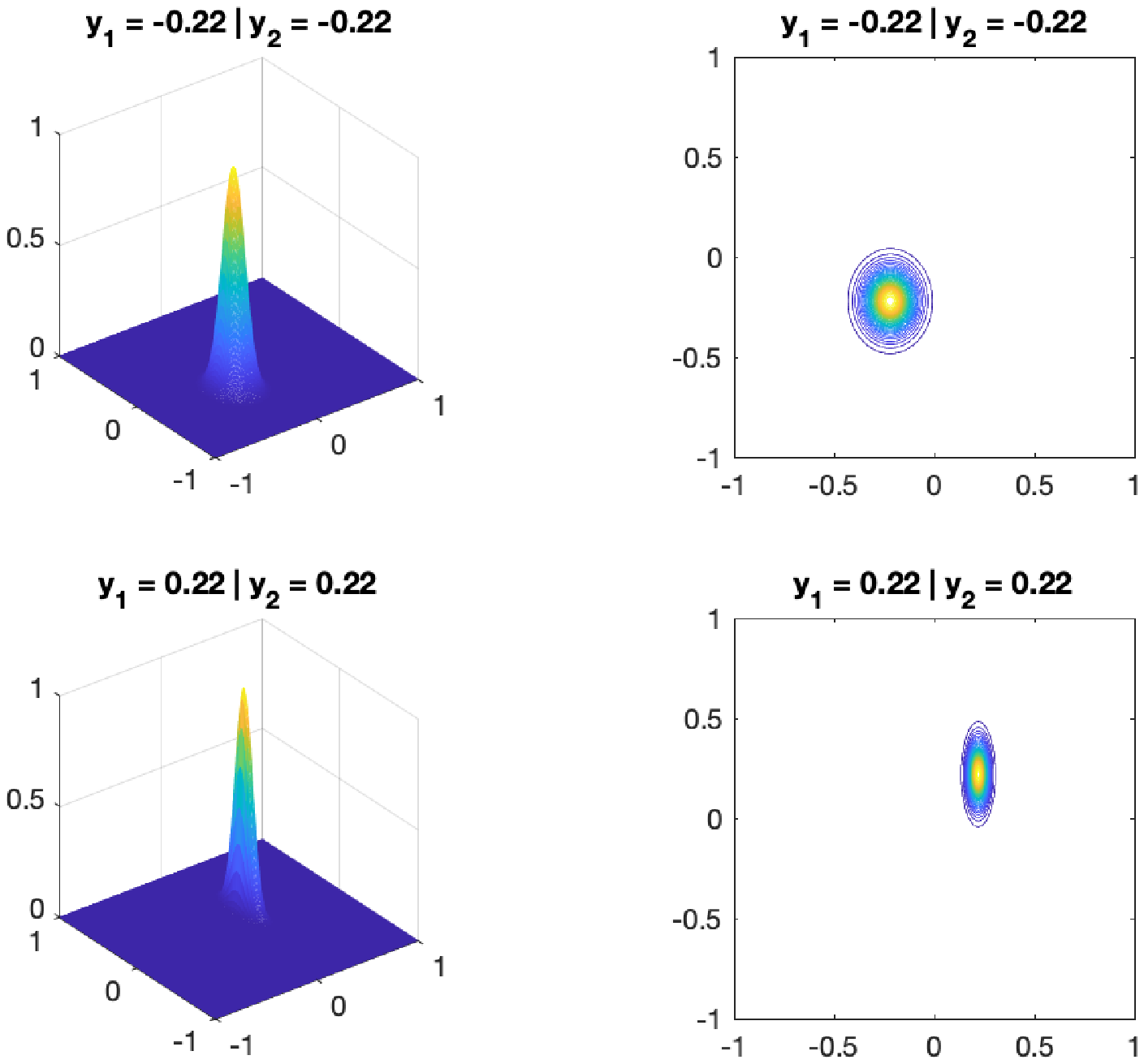}}
\caption{One peak problem solutions on the unit domain: $\alpha=1.54$ (top),  $\alpha=9.46$ (bottom).}
\label{fig:sc4.3sols}
\end{figure}

In the paper~\cite{LangSS20}, the one peak test problem defined on the unit domain 
is made {\it anisotropic} by  scaling the solution in  the first coordinate
direction  by a linear function $\alpha(y_1)= 18 y_1 + 11/2$ so that 
$\alpha$ takes values in the interval $[\rblx{1},10]$. The corresponding pathwise solution
is then given by 
\begin{align} \label{peaksol}
u(x,\y) &= \exp ( - \rblx{50}  \{  \alpha(y_1) (x_1 -y_1)^2 + (x_2 -y_2)^2 \} ).
\end{align}
\rblx{The solution \eqref{peaksol} is generated  by specifying an uncertain forcing function
\begin{subequations}
   \begin{align} \label{eq:dscaled}
      f(x,\y) &= d(x_1,x_2,y_1,y_2) \cdot  \exp ( - \beta  \{  \alpha(y_1) (x_1 -y_1)^2 + (x_2 -y_2)^2 \} )
   \\
   \intertext{with}
   \label{eq:fscaled}
      d(x_1,x_2,y_1,y_2) &= -4\beta^2 \left \{ \alpha^2(y_1)  (x_1 -y_1)^2 + (x_2 -y_2)^2 \right \}  
      + 2\beta (\alpha(y_1)+1) .
   \end{align}
\end{subequations}
Realisations of  the reference solution~\eqref{peaksol}   are shown
at two distinct sample  points  in Fig.~\ref{fig:sc4.3sols}.}
The anisotropy introduced by the \rblx{scaling with $\alpha$} is  a clear  feature.

Our specific  goal  is to compute the following quantity of interest (QoI)
\begin{align} \label{QoIref}
\Bbb{E}\left[\phi_{\rblx{I}}(u)\right] &=   \int_{\rblx{\left[-\frac 14, \frac 14\right]^2}} \int_I u^2(x,\y) \, \dx \, \rblx{\mathrm{d} \pi(\y)} , 
\end{align}
where $\phi_{\rblx{I}}(u)=  \int_I u^2(x,\cdot) \, \dx$.
The choice $\beta=50$ is then helpful for two reasons:
\begin{itemize}
\item 
The Dirichlet boundary condition ($u$ satisfying \eqref{peaksol} on \rblx{$\partial I$}) may be replaced
 without significant loss of accuracy by the numerical approximation 
$ u_{\bullet \z}=0$ on \rblx{$\partial I$}.
\item
A reference value  (accurate to more than 10 digits)
\begin{align} \label{QoIref:10:digits}
   \Bbb{E}\left[\phi_{\rblx{I}}(u)\right]\approx  Q := {1\over 9} \cdot (\sqrt{10 }-1)\cdot {\pi\over \beta} = 0.015 095 545 \ldots
\end{align}
 may be readily computed\rblx{; see~\cite[Appendix]{LangSS20}} for details. 
 \end{itemize}

\begin{figure}[!thp]
\centering
\includegraphics[width = 0.8\textwidth]{{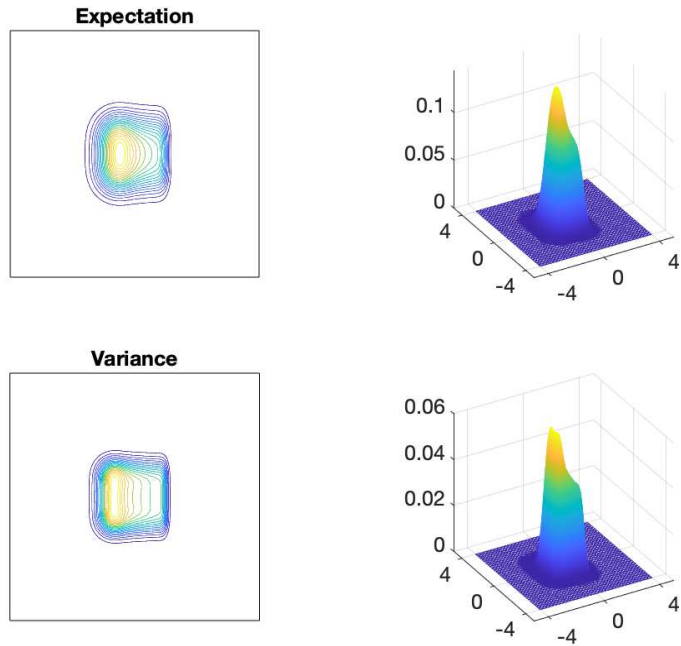}}
\caption{Reference solution for test case III.}
\label{fig:sc4.3refsol}
\end{figure}

\noindent
\rblx{We compute estimates of the QoI by solving the problem~\eqref{eq:pde2:strong}
using the coordinate transformations $x_j\leftarrow 4 x_j$ and $y_j\leftarrow 4 y_j$ ($j = 1,2$).
In this case, the pathwise solution  on the scaled domain
$D \times \Gamma$ is given by~\eqref{peaksol} by specifying $\beta = 50/16$
and $\alpha(y_1)= (9 y_1 + 11)/2$.
Moreover, the QoI in~\eqref{QoIref} (and its reference value given in~\eqref{QoIref:10:digits})}
can be estimated within Algorithm~\ref{algorithmx} by computing \rblx{the following~quantity:}
\begin{align*} 
   \rblx{\frac {1}{16}}\, \Bbb{E}\left[\phi_{\rblx{D}}(\rblx{u_{\ell}^{\rm SC}})\right] &=
   {1\over 16}\, \int_\Gamma \int_D \rblx{\big( u_{\ell}^{\rm SC}(x,\y) \big)^2}\, \dx \, \rblx{\mathrm{d} \pi(\y)} .
\end{align*}
A reference solution to the scaled problem is shown in Fig.~\ref{fig:sc4.3refsol}.

\begin{figure}[!thp]
\centering
\includegraphics[width = 0.8\textwidth]{{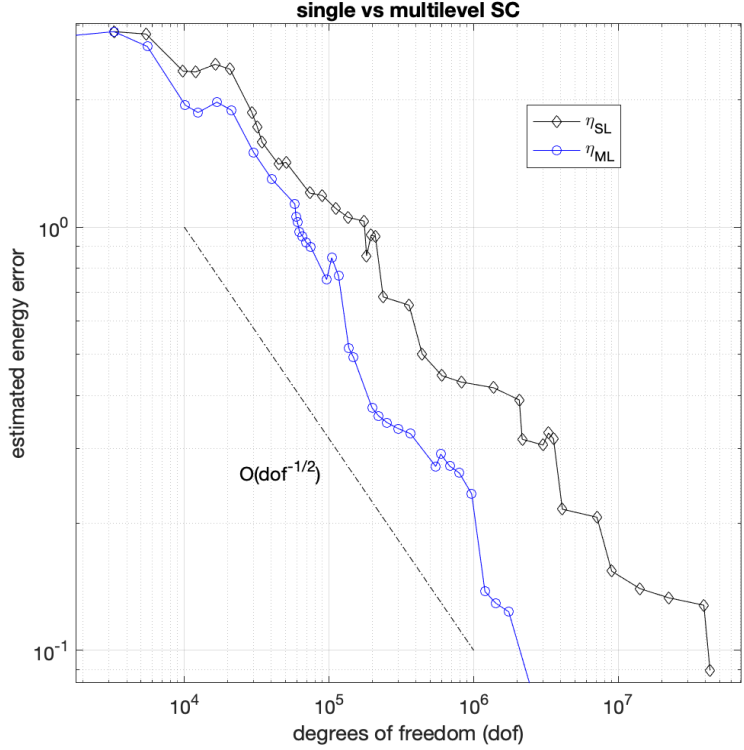}}
\caption{Evolution of the \rblx{single-level} and multilevel  error estimates  
for the one peak test problem with error tolerance set to {\tt 1e-1}.}
\label{fig:sc4.3errors}
\end{figure}

\begin{figure}[!pth]
\centering
\includegraphics[width = 0.32\textwidth]{{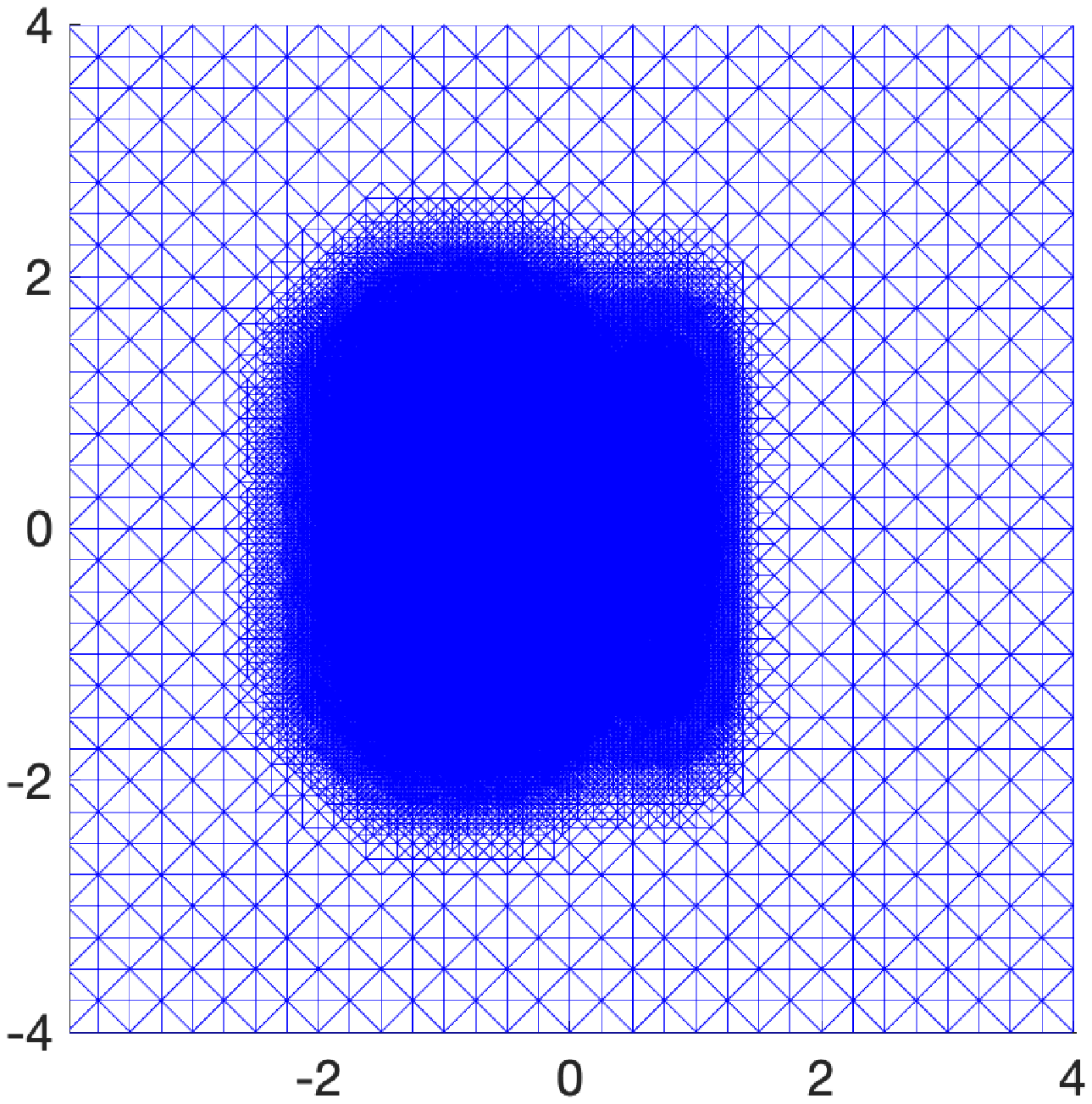}}
\includegraphics[width = 0.32\textwidth]{{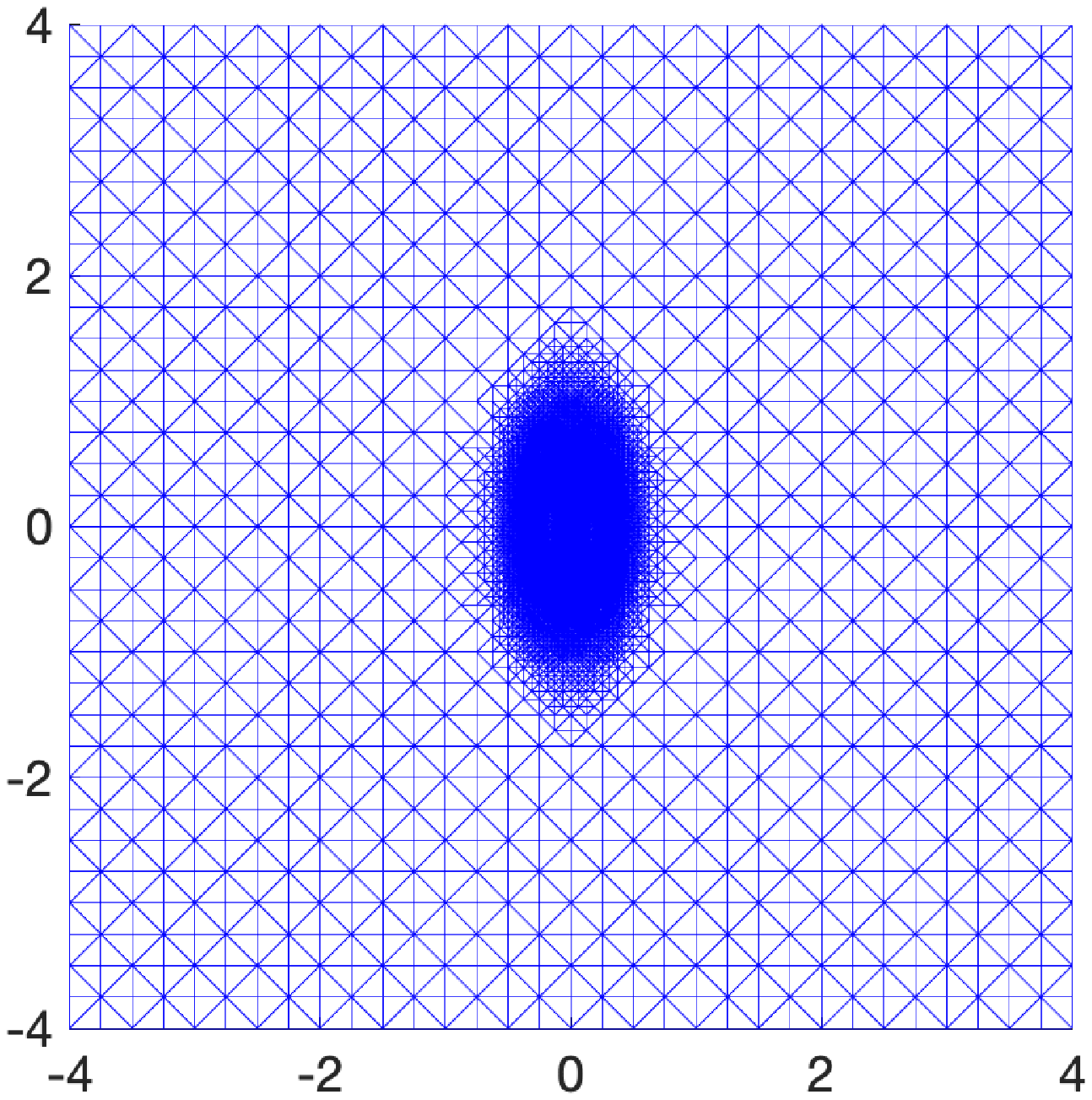}}
\includegraphics[width = 0.32\textwidth]{{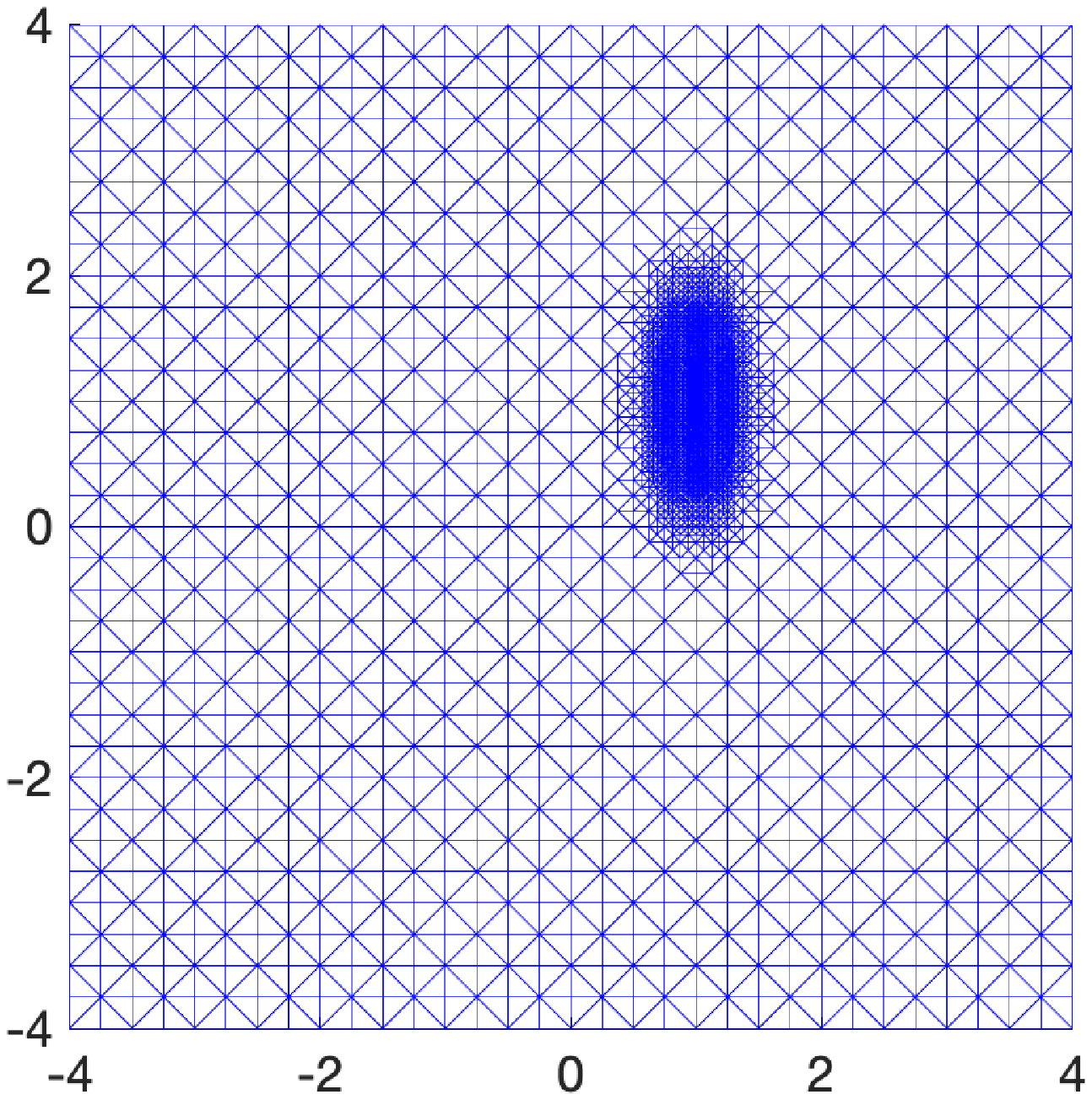}}
\caption{Single-level mesh  (left) and  meshes associated with the central 
collocation point (middle) and top right corner  point (right)  
when the tolerance is reached for test case III.}
\label{fig:sc4.3meshes}
\end{figure}

A comparison  of the  \rblx{single-level} and  multilevel \rblx{SC-FEM} algorithms when
applied to the one peak test problem is given by the \rblx{evolution of error estimates} in Fig.~\ref{fig:sc4.3errors}.
The single-level algorithm reached the tolerance in 37 steps with 169 active collocation points and the final 
approximation had {\tt 42,961,659} degrees of freedom. The multilevel algorithm proved to be
much more efficient. The same tolerance was reached  in 34 steps with 153 collocation points 
in the final approximation space\rblx{. Crucially,} each \rblx{collocation point is} associated with  a mesh
 that is locally refined 
in the vicinity of \rblx{the respective point in~$D$} (as illustrated in Fig.~\ref{fig:sc4.3meshes}).
\rblx{In contrast, the final mesh generated by the adaptive single-level SC-FEM}
has refinement everywhere in \rblx{a larger region corresponding to the union of supports of all sampled solutions}.
When the error tolerance \rblx{was} reached, both algorithms gave estimates of the QoI that agreed with the
reference value \rblx{to five} decimal places (0.015092 for the single-level case vs 0.015087 for the multilevel case).

The upshot of the effective use of tailored refinement is an order of magnitude decrease in the
overall computation time. The total number of degrees of freedom in the multilevel case
was \rblx{{\tt 2,620,343}---a} factor of \rblx{16} reduction overall.  Looking at the 
associated rates of convergence we see that the optimal rate $O({\rm dof}^{-1/2})$
is  recovered in the multilevel case. We anticipate that similar performance gains 
will be realized whenever a problem has local features that can be effectively 
resolved using  sample-dependent meshes.

We have also solved the one peak test problem using an efficient adaptive stochastic Galerkin approximation strategy.
While the linear algebra associated with the
Galerkin formulation is decoupled in \abrevx{this} case, the computational overhead of evaluating
the right-hand side vector is a significant limiting factor in terms of the relative efficiency.
The overall CPU time taken to compute 4 digits in the QoI using adaptive \abrevx{stochastic Galerkin FEM} is
comparable to the CPU  time taken to compute 5 digits using the multilevel \abrevx{SC-FEM}~strategy.

\section{Conclusions}\label{sec:conclusions}

Adaptive methods hold the key to efficient approximation of solutions to linear elliptic partial differential \rblx{equations} 
with random data. The numerical results presented in this series of two  papers demonstrate the effectiveness 
and the robustness of our novel SC-FEM error estimation strategy, as well as the utility of the error indicators guiding
 the  adaptive refinement process. Our results also suggest that optimal rates of convergence are more
difficult to achieve in a sparse grid collocation framework than in a 
multilevel stochastic Galerkin framework. It is demonstrated herein that
the overhead of generating  specially tailored sample-dependent meshes can be  worthwhile 
and optimal convergence rates can be recovered when the solutions to the sampled problems have 
local features in space. The single-level strategy discussed in part~I of this work is, however, likely to be  more 
efficient (certainly in terms of overall CPU time) when a single adaptively refined grid can 
adequately resolve  spatial features associated with solutions to a range of individually sampled problems.

\bibliographystyle{siam}
\bibliography{references}

\end{document}